\documentclass[10pt, a4paper]{article}
\usepackage[utf8]{inputenc}
\usepackage{microtype} 
\usepackage[english]{babel}
\usepackage{authblk}
\usepackage{amsmath}
\usepackage{amsthm}
\usepackage{amsfonts}
\usepackage{amssymb}
\usepackage{graphicx}
\usepackage{subcaption}
\usepackage{mathrsfs}
\usepackage{enumerate}
\usepackage{mathtools} 
\usepackage{xfrac} 
\usepackage{geometry}
\usepackage{appendix}
\usepackage{authblk}
\usepackage{booktabs}
\usepackage{multirow}
\usepackage{float} 
\usepackage{tikz} 

\usepackage[normalem]{ulem} 

\usepackage[autostyle]{csquotes} 
\usepackage[
  backend=biber,
  style=alphabetic,
  giveninits=true, 
  terseinits=true, 
  maxnames=6, 
  natbib=true, 
  sortlocale=en\_GB,
  url=false,
  doi=true,
  isbn=false,
  eprint=true
  ]{biblatex}
\usepackage{hyperref}
\usepackage{xcolor}
\hypersetup{
  colorlinks=true,
  linkcolor={blue!80!black},
  citecolor=black,
  urlcolor={blue!80!black}
}

\usepackage{algorithm} 
\usepackage{algorithmic} 

\theoremstyle{plain}
\newtheorem{theorem}{Theorem}[section]

\theoremstyle{definition}

\theoremstyle{remark}

\numberwithin{equation}{section}

\newcommand{\R}{\mathbb{R}}

\newcommand{\supp}{\operatorname{supp}}

\DeclarePairedDelimiter{\abs}{\lvert}{\rvert}

\newcommand{\sga}{a} 
\newcommand{\klor}{r} 
\newcommand{\kloT}{T} 

\title{Statistical comparison of reconstruction methods for the inverse boundary problem of the one-dimensional wave equation} 

\author{Samuel Agenorwoth}
\affil{Department of Computational Engineering, LUT University, Lahti, Finland}

\author{Emilia L.K. Bl{\aa}sten\footnote{Corresponding author, \texttt{emilia.blasten@iki.fi}}}
\affil{Department of Computational Engineering, LUT University, Lahti, Finland}

\begin{document}

\maketitle

\abstract{
  Several numerical reconstruction algorithms for the inverse boundary value problem of the 1-dimensional wave equation exist. In this paper we revisit two of them, the Sondhi--Gopinath (SG) method from 1971 and the Korpela--Lassas--Oksanen (KLO) method from 2016. Our main contribution is a systematic and statistical comparison of the performance of these two well established methods. The stability of the SG algorithm, even in practical applications, is well known. KLO has a regularisation scheme with a theoretical proof, and is an evolution of the boundary control method. Both are based on the idea of constructing solutions that are characteristic functions of a set at a given time. This similarity has been pointed out before, but no systematic comparison has been published.

  We compare the performance of the two algorithms with noisy simulated data. The application in our mind is the reconstruction of the internal cross-sectional area of a pressurised fluid pipe. The area corresponds to  the first order $\partial_x$-term of the wave equation. Instead of just observing the performance on a few test cases, we generate $n=1000$ random area profiles of various smoothness levels and measurement noise up to $10\%$ of the signal energy and perform statistical tests. SG and KLO have a difference of one time-derivative in their standard boundary data, which complicates the analysis. Our results show that SG performs better in the low noise regime, and KLO with high noise. SG is easier to implement and runs faster.
}

\bigskip {\noindent \bf Keywords:} inverse problems, numerical comparison, noisy data, wave equation


\section{Introduction}
The inverse problem for the first-order term of the one-dimensional wave equation involves solving the unknown coefficient function $A(x)$ from the measured $H(0,t)$ from either of
\begin{equation}
  \label{eq:intro-wave-eq}
  \begin{cases}
    \partial_t H(x,t) = - \frac{1}{A(x)} \partial_x Q(x,t) \\
    \partial_t Q(x,t) = -A(x) \partial_x H(x,t) \\
    H(x,0) = Q(x,0) = 0 \\
    Q(0,t) = \delta_0(t)
  \end{cases}
  \qquad \mbox{or} \qquad
  \begin{cases}
    \partial_t^2 H(x,t) = \partial_x^2 H(x,t) + \frac{\partial A(x)}{A(x)} \partial_x H(x,t) \\
    H(x,0) = \partial_t H(x,0) = 0 \\
    A(0) \partial_x H(0,t) = f(t)
  \end{cases}
\end{equation}
where $f$ can be varied freely. Due to applications involving pressurised water pipes \cite{ghidaoui05_review_water_hammer_theor_pract} or acoustic wave guides \cite{titze00_princ_voice_produc}, we call $A(x)$ the cross-sectional area at points $x \in (0,\infty)$, $H$ the fluid pressure, and $Q$ the flow rate.

The problem has been studied and solved with various methods in the last 70 years \cite{gelfand55_deter_differ_equat_from_its_spect_funct,sondhi71_deter_vocal_shape_from_impul,korpela16_regul_strat_inver_probl_1}. Past results range from spectral to dynamic reflection data to translating between them \cite{gelfand55_deter_differ_equat_from_its_spect_funct,sondhi71_deter_vocal_shape_from_impul,blagoveshchenskii71_local_method_solut_nonst_inver,browning99_time_frequen_domain_scatt_one,browning00_time_frequen_domain_scatt_one,belishev08_bound_contr_inver_probl}, to transmission problems \cite{rakesh96_imped_inver_from_trans_data_wave_equat,rakesh00_charac_trans_data_webst_horn_equat}, to numerical reconstruction and regularisation strategies \cite{korpela16_regul_strat_inver_probl_1}, to range characterisation \cite{belishev01_dynam_system_with_bound_contr,belishev16_charac_inver_data_bound_contr_method}, and applications to real-world noisy data \cite{sondhi83_inver_probl_vocal_tract}. One could claim that this setting is the most well-understood non-linear inverse problem! Many of the methods spawned during its research have led to generalisations across higher dimensions, e.g. the BC method \cite{belishev87_approac_to_multid_inver_probl_wave_equat,belishev97_bound_contr_recon_manif_metric_bc} and variants \cite{oksanen11_solvin_inver_probl_wave_equat,korpela16_regul_strat_inver_probl_1}.

Our contributions to the understanding of this classical problem are the following. We made a statistical comparison of two algorithms, one from \cite{sondhi71_deter_vocal_shape_from_impul} and the other from \cite{korpela16_regul_strat_inver_probl_1}, comparing their performance under various noise levels and smoothness assumptions. The main novelty is first of all that this comparison even happens, but also that our method of comparison makes sure that the algorithms are tested on a wide range of unknown functions and noise levels, and a hypothesis test confirms which one is better. Furthermore, to do the above fairly, we have implemented both algorithms only as described in their papers, without contact from their authors. As a consequence we have reproducibility of the much more recent algorithm of \cite{korpela16_regul_strat_inver_probl_1}. The results say that KLO is better in the high noise regime, but SG runs faster, is easier to implement, and runs better in the low noise regime. Smoothness levels of the unknown function had no appreciable effect one way or the other in this comparison.

  We will continue with a motivation for this research. Then in the following sections we start with a comparison of the theory of the two methods, then set up our Monte-Carlo regime, continue with numerical implementation details, and finally, present the results.

\medskip
The equations (\ref{eq:intro-wave-eq}) can be used to model situations which are inherently one-dimensional such as sound propagation in the vocal tract \cite[chapter 6]{titze00_princ_voice_produc,lukkari12_webst_equat_with_curvat_dissip} or fluid pressure in water supply pipes \cite{ghidaoui05_review_water_hammer_theor_pract}, or situations where a higher-dimensional model reduces to a one-dimensional equation \cite{rakesh98_inver_spher_symmet_poten_from}. Despite the extensive theoretical understanding of the one-dimensional hyperbolic inverse problem, these engineering and medical applications are far from being solved in practice. Novel cheaper measurement methods are being explored in engineering \cite{zouari19_inter_pipe_area_recon_as,blåsten19_block_detec_networ,zouari20_exper_verif_accur_robus_area,gozum22_noise_based_high_resol_time}, and the problem of inverse glottal filtering \cite{alku11_glott_inver_filter_analy_human,drugman14_glott_sourc_proces,alku19_openg_open_envir_evaluat_glott_inver_filter} is still open despite recent theoretical steps for determining the source and shape by one measurement \cite{feizmohammadi25_recon_1d_evolut_equat_their}.

What is unclear is how the existing one-dimensional methods perform with respect to each other. In this paper we will compare the methods introduced by \textcite{sondhi71_deter_vocal_shape_from_impul} and \textcite{korpela16_regul_strat_inver_probl_1} numerically (SG and KLO for short). The latter method is an evolution of the boundary control (BC) method \cite{belishev97_bound_contr_recon_manif_metric_bc,bingham08_iterat_time_rever_contr_inver_probl,oksanen11_solvin_inver_probl_wave_equat}. According to \cite{belishev97_bound_contr_recon_manif_metric_bc}, Rakesh pointed out the similarities between the Sondhi--Gopinath method and one-dimensional BC method \cite[preprint in 1996]{belishev08_bound_contr_inver_probl}. How similar are they? What about computational performance, or stability to noise in the measurements? To our knowledge, no systematic comparison between these two methods have been published.

\smallskip We compare both methods by implementing them numerically and doing a Monte Carlo simulation by randomly generating the unknown coefficient function with various levels of smoothness  \cite{williams2006gaussian} and additive Gaussian measurement noise. In contrast to most studies with numerical implementations of reconstructions algorithms for inverse boundary value problems, we don't just observe the quality of the solution in a few test cases, but perform general statistical tests. As a result our conclusions are statistically significant in the categories of coefficient function under consideration.

In both SG and KLO the boundary pressure $H(0,t)$ is measured. Rewriting the 1st order system in \eqref{eq:intro-wave-eq} as a 2nd order scalar equation, the two direct problems are
\begin{equation}
  \label{eq:problem-intro}
  \begin{cases}
    \partial_t^2 H = \partial_x^2 H + \frac{\partial A(x)}{A(x)} \partial_x H, & \\
    H=0, & t \leq 0 \\
    -A(0) \partial_x H(0,t) = \partial_t \delta_0(t), & t\in \R
  \end{cases}
  \qquad \mbox{ and }\qquad
  \begin{cases}
    \partial_t^2 H = \partial_x^2 H + \frac{\partial A(x)}{A(x)} \partial_x H, & \\
    H=0, & t \leq 0 \\
    A(0) \partial_x H(0,t) = f(t), & t\in \R    
  \end{cases}
\end{equation}
where $f$ is an input freely controllable. The standard data of SG (on the left) is the impulse response function $\delta_0(t) \mapsto H(0,t)$, while KLO (on the right) has the Neumann-to-Dirichlet operator $f(t) \mapsto H(0,t)$. Due to the one-dimensional geometry, this operator can be determined by simply knowing the value of $H(0,t)$ when $f$ is the delta distribution $\delta_0(t)$. One notes an extra time-derivative in the SG problem on the left. Both of these are natural measurements. For example setting $Q(0,t)=\delta_0(t)$ in SG corresponds to creating a strong explosive sound, as described in \textcite{sondhi83_inver_probl_vocal_tract}, but setting $A(0)\partial_xH(0,t)=\delta_0(t)$ as in KLO (equivalent to $\partial_tQ(0,t) = - \delta_0(t) / A(0)$) can be done by opening a valve instantly (e.g. by hitting it with a sledgehammer) and leaving it open as in \cite{zouari20_exper_verif_accur_robus_area}. Should one wish to use SG to solve a problem with a KLO-type measurement, one would need to differentiate the measured $H(0,t)$ before applying SG. This leads to disaster if not done carefully when the measurement is noisy.

\smallskip The results of our study show a clear noise-dependent distinction between the two reconstruction strategies. In noise-free synthetic experiments, SG gives the most accurate reconstructions, particularly for smooth and moderately rough profiles, and consistently outperforms KLO in both $L^2$- and $H^1$-based error metrics. However, under additive measurement noise, the statistical comparison reveals a crossover: KLO becomes more robust than SG at moderate and high noise levels, especially in the derivative sensitive $H^1$-metric. Both methods are challenged by highly irregular hybrid profiles with sharp local transitions, although the dominant geometric structure is still recovered.

From a computational and implementation perspective, SG inverse solver is the simpler and more efficient method. Its Fredholm-type reconstruction can be implemented using Toeplitz or Woodbury acceleration, leading to relatively low memory requirements. KLO, in contrast, requires construction of the Neumann-to-Dirichlet map, assembly or application of the connecting operator, and repeated regularised projected solves over the reconstruction radius, making it more expensive and more sensitive to implementation details.

Based on our experiments, SG is recommended when the data are clean, the geometry is sufficiently regular, and computational simplicity is important. KLO is preferable when the available boundary data are noisy and robustness to perturbations is more important than achieving the smallest possible noiseless error. Thus, the practical choice depends on the data regime: SG is better high-accuracy method for ideal or weakly perturbed synthetic data, while KLO provides the more stable reconstruction strategy for contaminated measurements provided that the inlet regularisation and grid-resolution effects are carefully controlled.

\smallskip The structure of the paper is as follows: In Section~\ref{section2:comparison} we recall the two reconstruction algorithms and describe the methods we use for generating random area functions. In Section~\ref{section3:implementation} we describe the numerical implementation and comparison experiments. Sections~\ref{section4:results} and \ref{section5:conclusion} summarise the results and conclusions.

\section{Comparing two reconstruction algorithms}
\label{section2:comparison}

In this section, we start by describing and comparing the theoretical aspects of SG and KLO. We then describe the methodology used for generating area functions of a given smoothness randomly.

In the area reconstruction study, we generate three families of random area profiles: (i) smooth geometries based on a squared-exponential kernel, (ii) rough geometries based on Mat\'ern kernels with varying smoothness parameters, and (iii) hybrid geometries that  combine smooth and rough components with localized rectangular bumps. By varying the kernel choice and its parameters, we can systematically control the geometric complexity and assess its impact on reconstruction accuracy. We guarantee the strict positivity of the area by modelling its algorithm as a Gaussian process \cite{williams2006gaussian}.

The study of multiple variants of area functions with these different geometric characteristics allow us to assess how the smoothness, roughness, and presence of localized features affect reconstruction accuracy and the reliability of uncertainty estimates.

By running the reconstruction algorithm on many such profiles, we quantify the statistical distribution of reconstruction errors, assess the influence of geometric features i.e., smoothness, roughness, and local bumps, and validate the uncertainty estimates. The three families of generators are chosen to cover a wide range of geometric complexity, ensuring that our conclusions about the two reconstruction algorithm performance are not limited to a single, overly simplistic class of shapes.

\subsection{Reconstruction algorithms}

\label{section:compare-theory}

We will consider the problems \eqref{eq:intro-wave-eq} or equivalently \eqref{eq:problem-intro} for the wave equation on a one-dimensional segment $(0,\ell)$. An impulse of the form $Q(x,0)=\delta_0(t)$ (SG) or $A(0)\partial_xH(0,t)=f(t)$ (KLO) propagates into the fluid from an inlet on the left $x=0$ and the equations keep track of the pressure $H(x,t)$ and flow velocity $Q(x,t)$. We assume the wave speed is constant $c=1$, and consider varying cross-sectional area $A(x)$ that's constant when $x \geq \ell$.

In both cases we measure the pressure at the inlet, $H(0,t)$. We define the measurement operators
\begin{equation}
  \label{eq:measurement}
  \Lambda_{SG}\delta_0(t) = H(0,t), \qquad \Lambda_{KLO}f(t) = H(0,t), \qquad t \in \R
\end{equation}
where $H$ solves the left side of \eqref{eq:intro-wave-eq} or \eqref{eq:problem-intro} for SG and the right side for KLO. The inverse problem is about recovering the unknown area function $A(x)$, $x \in (0,\ell)$ from the measurement operators. Note that these measurement contain the initial impulse,
\begin{equation}
  \label{eq:lambda-response}
  \Lambda_{SG}\delta_0(t) = \frac{1}{A(0)} \big( \delta_0(t) + h_{SG}(t) \big), \qquad \Lambda_{KLO} \delta_0(t) = \frac{1}{A(0)} \big( \chi_{(0,\infty)}(t) + h_{{KLO}}(t) \big)
\end{equation}
so we write $h_{SG}(t)$ or $h_{KLO}(t)$for the corresponding response.

Theoretically, the SG and KLO reconstruction algorithms have the same main point: given boundary measurements, compute boundary values which would create a wave that is the characteristic function of some set at a given time. This is done by relying only on boundary measurements, i.e. without knowing the underlying area function $A(x)$. In both methods, the computation means solving boundary integral equations multiple times. The algorithms are different though, and not just due to a time-derivative in the boundary data. For example KLO originates from higher dimensions, where one needs to use deep unique continuation results by Tataru \cite{tataru95_unique_contin_solut_to_pdes,tataru99_unique_contin_operat_partial_analy_coeff} instead of simply d'Alembert's formulas after flipping time and space $\partial_t^2-\partial_x^2 \leftrightarrow \partial_x^2-\partial_t^2$ which makes unique continuation work in one dimensions. Moreover KLO has a proven regularisation scheme \cite{korpela16_regul_strat_inver_probl_1}, while SG seems stable both with numerical noise and real world data \cite{sondhi83_inver_probl_vocal_tract,zouari19_inter_pipe_area_recon_as,blåsten19_block_detec_networ,zouari20_exper_verif_accur_robus_area}.

In the rest of this section, we will recall the algorithms after which we will go into the details of the numerical comparison.

\subsubsection{Sondhi--Gopinath 1971}
Given the impulse-response \eqref{eq:lambda-response} of the first order system on the left of \eqref{eq:intro-wave-eq}, the Sondhi--Gopinath algorithm~\cite{sondhi71_deter_vocal_shape_from_impul} solves for the unknown function $A$ in a range $(0,\ell)$ in the following way:
\begin{enumerate}
\item given any positive $\sga < \ell$,
\item we construct a boundary function $f_{\sga}$ based on $\Lambda_A$ such that if in \eqref{eq:intro-wave-eq} we replace the boundary condition by $Q(0,t) = f_{\sga}(t)$, then the solution $H$ satisfies
  \begin{equation}
    H(x,\sga) =
    \begin{cases}
      1, &\mbox{for $0 < x < \sga$},\\
      0, &\mbox{for $x > \sga$},
    \end{cases}
  \end{equation}
  at time $t=\sga$.
\item we recover the area at $x=\sga$ by
  \begin{equation}
    \label{eq:area-from-f-squared}
    A(\sga) = \left( f_{\sga}(\sga) \right)^2.
  \end{equation}
\item repeat from step 1 until $A(x)$ is determined on a dense enough grid.
\end{enumerate}

The determination of $f_{\sga}$ comes is based on a form of unique continuation in the theorem below. It is a direct consequence of d'Alembert's formulas in one dimension since we can simply change the role of $x$ and $t$.
\begin{theorem}
  If $\widetilde{H}$ and $\widetilde{Q}$ satisfy
  \begin{align*}
    \partial_t\widetilde{H}(x,t) &= - \frac{1}{A(x)} \partial_x\widetilde{Q}(x,t), \qquad t > 0,\, x \in (0,\infty) \\
    \partial_t\widetilde{Q}(x,t) &= - A(x) \partial_x \widetilde{H}(x,t), \qquad t > 0,\, x \in (0,\infty)
  \end{align*}
  and $\widetilde{H}(0,t) = 2$, $\widetilde{Q}(0,t)=0$ on $0 < t < 2\sga$, then $\widetilde{H}(x,t)=2$ and $\widetilde{Q}(x,t)=0$ in the space-time triangle $0< x < \sga - \abs{\sga-t}$.
\end{theorem}
To find $f_{\sga}$, we use the ansatz
\begin{align}
  \widetilde{H}(x,t) &= H(x,t) + H(x,2\sga - t) \\
  \widetilde{Q}(x,t) &= Q(x,t) - Q(x,2\sga - t)
\end{align}
where $H,Q$ satisfy the 1st order system with initial conditions. Setting $x=0$, $Q(0,t)=f_{\sga}(t)$, using $\widetilde{Q}(0,t)=0$ and noting that
\begin{equation}
  H(0,t) = \Lambda_{SG}[Q(0,\cdot)](t) = \int \Lambda_{SG}[\delta_0](t-s) Q(0,s) ds,\label{convolution_of_SG}
\end{equation}
we arrive at the inhomogeneous Fredholm equation of the second kind
\begin{equation}
  \label{eq:fredholm}
  f_{\sga}(t) + \frac{1}{2} \int_0^{2\sga} h_{SG}(\abs{t-s}) f_{\sga}(s) ds = A(0), \qquad 0 < t < 2\sga
\end{equation}
for $f_{\sga}$

To reconstruct the area, we multiply the $\partial_tH$-equation in \eqref{eq:intro-wave-eq} by $A$, integrate over the space-time rectangle $(0,\sga) \times (0,\ell)$ and simplify using the initial conditions to get
\begin{equation}
  \int_0^{\sga} Q(0,t) dt = \int_0^{\ell} A(x) H(x,\sga) dx
\end{equation}
which leads to
\begin{equation}
  \label{eq:SG-volume}
  \int_0^{\sga} f_{\sga}(t) dt = \int_0^{\sga} A(x) dx.
\end{equation}
Further algebra applied to (\ref{eq:fredholm}) and (\ref{eq:SG-volume}) leads to (\ref{eq:area-from-f-squared}), see equations (12)--(16) in \cite{sondhi71_deter_vocal_shape_from_impul} for details.

\subsubsection{BC-algorithm from Korpela--Lassas--Oksanen 2016}

In \cite{oksanen11_solvin_inver_probl_wave_equat} the author considers an inverse problem for an $(n+1)$-dimensional wave equation analogous to
\begin{equation}
  \partial_t^2 u = \frac{1}{\mu(x) \sqrt{g(x)}} \partial \left( \frac{\mu(x)}{\sqrt{g(z)}} \partial_x u \right)
\end{equation}
and showed a reconstruction method based on boundary control. The method was regularised and numerically implemented later in $1$ spacial dimension in \cite{korpela16_regul_strat_inver_probl_1}. Setting $g(x)=1$, $\mu(x) = A(x)$ and $u=H$ gives our wave equation on the right side of \eqref{eq:intro-wave-eq}. In contrast to the Sondhi--Gopinath-type measurements of $\Lambda_{SG}[Q(0,\cdot)](t)=H(0,t)$, Korpela--Lassas--Oksanen defines
\begin{equation}
  \label{eq:KLO-lambda}
  A(0) \partial_x H(0,t) = f(t) \implies \Lambda_{KLO}f(t) = H(0,t).
\end{equation}
Comparing to $\Lambda_{SG}$, if we set $A(0) \partial_x H(0,t) = f(t)$, then
\[
  Q(0,t) = \int_0^t (\partial_t Q)(0,s) ds = - \int_0^t A(0) \partial_x H(0,s) ds = -\int_0^t f(s) ds,
\]
so
\begin{equation}
  \label{eq:KLO-SG-meas}
  \Lambda_{KLO}f = \Lambda_{SG} \left(-\int f \right) = - \int \Lambda_{SG} f
\end{equation}
and vice versa
\begin{equation}
  \label{eq:SD-KLO-meas}
  \Lambda_{SG} f = \Lambda_{KLO}(-\partial f) = -\partial_t \Lambda_{KLO} f.
\end{equation}

Given $\Lambda_{KLO}$, the Korpela--Lassas--Oksanen algorithm recovers the unknown function $A$ in the range $(0,\ell)$ as follows. Let $\kloT$ be a large enough time so that a wave from $(x,t)=(0,0)$ reaches $x=\ell$ in time $\kloT$. With our assumption of wave speed 1, we have $\kloT=\ell$. The algorithm:
\begin{enumerate}
\item given any positive $\klor < \kloT$,
\item \label{KLO-step2} we select $\alpha>0$ and find a boundary function $f_{\klor,\alpha}$ supported on $[\kloT-\klor,\kloT]$ such that if in the right side of \eqref{eq:intro-wave-eq} we replace the boundary condition by $A(0) \partial_x H^{f_{\klor,\alpha}}(0,t) = f_{\klor,\alpha}(t)$, then the energy
  \begin{equation}
    \label{eq:KLO-energy}
    \int_0^{\ell} \left| H^{f_{\klor,\alpha}}(x,\kloT) - 1 \right|^2 A(x) dx + \alpha \int_0^{2\kloT} \left| f_{\klor,\alpha}(t) \right|^2 dt
  \end{equation}
  is minimised. This can be done using only knowledge of $\Lambda_{KLO}$.
\item \label{KLO-step3} we recover the area at $x=\klor$ by
  \begin{equation}
    \label{eq:KLO-area}
    A(\klor) = \partial_{\klor} \lim_{\alpha \to 0} \int_{T-\klor}^T f_{\klor,\alpha}(t) (\kloT - t) dt
  \end{equation}
\item repeat from step 1 until $A(x)$ is determined on a dense enough grid.
\end{enumerate}

Steps 2 and 3 can be done without knowledge of the unknown area function $A(x)$, $x\in (0,\ell)$ \cite{korpela16_regul_strat_inver_probl_1,oksanen11_solvin_inver_probl_wave_equat}. This relies on the Blagoveshchenskii identities \cite{blagoveshchenskii69_one_dimen_inver_bound_value,blagoveshchenskii71_inver_bound_value_probl_theor}, written as in \cite{korpela16_regul_strat_inver_probl_1}, where $B1$ is the operator $B$ from \eqref{def-B-operator} applied to the constant function $1$. In fact $B1(t) = \chi_{(0,T)}(t)(T-t)$.
\begin{align}
  \int_0^{\ell} H^f(x,\ell) A(x) dx &= \int_0^{2\kloT} f(t) B1(t) dt \label{blago1}\\
  \int_0^{\ell} H^f(x,\ell) H^g(x,\ell) A(x) dx &= \int_0^{2\kloT} f(t) Kg(t) dt \label{blago2}
\end{align}
where $H=H^f$ satisfies the right side of \eqref{eq:intro-wave-eq}, similarly for $H^g$ with $f$ replaced by $g$ in the boundary condition. The operators involved are given by
\begin{align}
  Bf(t) &= \chi_{(0,\kloT)}(t) \int_t^{\kloT} f(s) ds, \label{def-B-operator}\\
  P_{\klor}f(t) &= \chi_{(\kloT-\klor,\kloT)}(t) f(t),\\
  Jf(t) &= \frac{1}{2} \int_0^{2\kloT} \chi_{\Delta}(t,s) f(s) ds = \frac{1}{2} \chi_{(0,\kloT)}(t) \int_t^{2\kloT -t} f(s) ds,\\
  Rf(t) &= f(2\kloT - t),\\
  Kf &= (\Lambda_{KLO}^{\ast}J - J \Lambda_{KLO})f = (R \Lambda_{KLO} R J - J \Lambda_{KLO})f,\label{Kf}
\end{align}
where $\Lambda_{KLO}^{\ast}$ is the $L^2(0,2\kloT)$-dual of $\Lambda_{KLO}$. Since $\abs{H^f-1}^2 = H^f H^f - 2 H^f \cdot 1 + 1$, minimising \eqref{eq:KLO-energy} is equivalent to minimising
\begin{align}
  &\int_0^{\ell} H^{f_{\klor,\alpha}}(x,\kloT) H^{f_{\klor,\alpha}}(x,\kloT) A(x) dx - 2 \int_0^{\ell} H^{f_{\klor,\alpha}}(x,\kloT) A(x) dx + \alpha \int_0^{2\kloT} \abs{f_{\klor,\alpha}(t)}^2 dt \notag \\
  &\qquad = \int_0^{2\kloT} f_{\klor,\alpha}(t) K f_{\klor,\alpha}(t) dt - 2 \int_0^{\kloT} f_{\klor,\alpha}(t) B1(t) dt + \alpha \int_0^{2\kloT} \abs{f_{\klor,\alpha}(t)}^2 dt \label{eq:KLO-energy-bilin-form}
\end{align}
By looking at the Fréchet derivative of \eqref{eq:KLO-energy-bilin-form}, according to \cite{korpela16_regul_strat_inver_probl_1,oksanen11_solvin_inver_probl_wave_equat}, this is minimised by
\begin{equation}
  \label{eq:KLO-minimiser}
  f_{\alpha,\klor} = (P_{\klor} K P_{\klor} + \alpha)^{-1} P_{\klor}B1
\end{equation}
taking into account the restriction $\supp f_{\klor,\alpha} \subset [\kloT-\klor,\kloT]$. Each of the operators involved here can be computed given the boundary measurements $\Lambda_{KLO}$. This allows computing $f_{\klor,\alpha}$ in Step~\ref{KLO-step2}.

Lastly, Step~\ref{KLO-step3} follows from a few facts: Firstly recall that $\supp f_{\klor,\alpha} \subset [\kloT-\klor,\kloT]$, so the constant wave speed of 1 implies that at time $t=\kloT$ we have
\begin{equation}
  \label{eq:KLO-H-support}
  \supp H^{f_{\klor,\alpha}}(\cdot, \kloT) \subset [0,\klor].
\end{equation}
Secondly,
\begin{equation}
  \label{eq:KLO-minimiser-almost-one}
  \lim_{\alpha\to0} \int_0^{2\ell} \abs*{H^{f_{\klor,\alpha}}(x,\kloT) -1}^2 A(x) dx = 0
\end{equation}
as seen in \cite{korpela16_regul_strat_inver_probl_1,oksanen11_solvin_inver_probl_wave_equat}. These two imply that $H^{f_{\klor,\alpha}}(\cdot,\kloT) \approx \chi_{(0,\klor)}$ in the $L^2(0,\ell;A)$-sense. This fact combined with the definition of the operator $B$ and \eqref{blago1} gives
\begin{align}
  \int_{\kloT-\klor}^{\kloT} f_{\klor,\alpha}(t) (\kloT - t) dt
  &= \int_0^{\kloT} f_{\klor,\alpha}(t) (\kloT - t) dt = \int_0^{2\kloT} f_{\klor,\alpha}(t) B1(t) dt = \int_0^{\ell} H^{f_{\klor,\alpha}}(x,\kloT) A(x) dx \notag \\
  &\to  \int_0^{\ell} \chi_{(0,\klor)}(x) \cdot A(x) dx =  \int_0^{\klor}A(x) dx 
\end{align}
as $\alpha \to 0$. Differentiation would give $A(\klor)$, however measurement noise might prevent the use of $\alpha=0$ in practical computations. Instead of $\alpha \to 0$ in \eqref{eq:KLO-area}, Korpela--Lassas--Oksanen \cite{korpela16_regul_strat_inver_probl_1} suggest using the regularisation strategy $\alpha = 2^{13/9} \cdot \ell^{4/9} \varepsilon^{4/9}$ where $\varepsilon$ is an upper bound for the noise between the theoretical $\Lambda_{KLO}$ and the observed one measured in operator norm of bounded operators acting on $L^2(0,\ell)$. We will explore this in sections \ref{Noise_4_KLO} and \ref{Noise_4_KLO_results}.

\subsection{Random area profile generation}\label{area_formulation}
Let the grid be \(x_0,x_1,\cdots,x_N\) with \(x_0=0\) and \(x_N=\ell.\) We model \(g(x) = \ln A(x)\) as a zero-mean Gaussian process, 
\begin{equation}
\mathbf{g} \sim \mathcal{N}(\mathbf{0},\mathbb{K}),\quad \mathbb{K}_{ij} = k(|x_i -x_j|),
\end{equation}
with a small jitter term $\varepsilon = 10^{-10}$ for numerical stability.  The area profile is \(A(x) = A_0\exp(g(x)),\) normalised so that $A(0)=1$ and clipped to $[A_{\min},A_{\max}]$. Here $\mathcal{N}(\boldsymbol{\mu},\boldsymbol{\Sigma})$ denotes
the multivariate normal distribution with mean $\boldsymbol{\mu}$ and covariance
$\boldsymbol{\Sigma}$, and $\mathbf{g} = [g(x_0),\dots,g(x_N)]^{\top}$ collects
the nodal values of the process. The function $k$ is the covariance function, or kernel, and $\mathbb{K}$ is the corresponding covariance matrix. In the kernels
below, $\sigma$ is the marginal standard deviation of the log-area field and $\overline\ell$ is its correlation length, so that values separated by much more than $\overline\ell$ are close to uncorrelated.

\subsubsection{Generative models}
\label{Covariance models}
We consider three generative models for the log-area field, each designed to explore different aspects of the reconstruction algorithm's behavior under uncertainty. The first two utilise standard covariance kernels within the Gaussian process framework while the third is composite. In what follows, we denote the distance between two spatial points by \(\overline r=|x-x'|.\)

\begin{enumerate}
\item \textbf{Squared-exponential (SE) kernel:}
\begin{equation}
k_{\text{SE}}(\overline r)  = \sigma^2 \exp\left(-\frac{\overline r^2}{2 {\overline\ell}^2}\right),
\end{equation}
producing infinitely differentiable sample paths. 

\item \textbf{Mat\'ern kernel:}
\begin{equation}
k_{\nu}(\overline r) =
\sigma^2\,\frac{2^{1-\nu}}{\Gamma(\nu)}
\left(\frac{\sqrt{2\nu}\,\overline r}{\overline\ell}\right)^{\nu}
K_{\nu}\!\left(\frac{\sqrt{2\nu}\,\overline r}{\overline\ell}\right),
\end{equation}
with closed forms at $\nu\in\{\tfrac{1}{2},\tfrac{3}{2},\tfrac{5}{2}\}$ yielding continuous but non-differentiable ($\nu=\tfrac{1}{2}$), once differentiable ($\nu=\tfrac{3}{2}$), and twice differentiable ($\nu=\tfrac{5}{2}$) sample paths \cite{stein1999interpolation,williams2006gaussian}. Of the three closed forms above, the inversion experiments reported below use $\nu = 1.5$ as the representative Mat\'ern case. Realisations with
$\nu = 0.5$ are excluded because they are continuous but not differentiable at any point, so the $H^1$ norm appearing in the denominator of our reported errors
does not converge under grid refinement, and reported relative $H^1$ errors would then reflect the discretisation rather than the quality of the
reconstruction. The case $\nu = 2.5$ is intermediate in smoothness between the squared-exponential and $\nu = 1.5$ profiles, both of which are reported.
Profiles rougher than any Mat\'ern realisation are covered by the hybrid model, whose rectangular components introduce jump discontinuities.

\item \textbf{Hybrid model:}
\begin{equation}
g(x) = g_{s}(x) + g_{\overline r}(x) + g_{\mathrm{rect}}(x),
\end{equation}
combining a smooth SE component $g_s$, a rough Mat\'ern component $g_{\overline r}$, and deterministic rectangular bumps $g_{\mathrm{rect}}(x)=\sum_{m=1}^{n_{\mathrm{rect}}} h_m\,\mathbf{1}_{[c_m-w_m/2,\, c_m+w_m/2]}(x)$ with random centres, widths and heights.
\end{enumerate}

Table \ref{tab:params} summarises the parameters chosen for each area generative model. The computational domain is defined as $x\in[0,\ell]$, with length $\ell =2$ for this specific test case. Areas are clipped to $[0.5,2.0]$ after normalisation.
\begin{table}[H]
\centering
\begin{tabular}{|l|l|c|}
\hline
\multicolumn{1}{|c|}{Generator} & 
\multicolumn{1}{|c|}{Parameter} & 
\multicolumn{1}{|c|}{Value} \\
\hline
\multirow{2}{*}{SE}     & $\overline\ell$   & $0.12$ \\
\cline{2-3} 
        & $\sigma$ & $0.20$ \\
\hline
\multirow{2}{*}{Mat\'ern ($\nu\in\left\{0.5,1.5,2.5\right\}$)}  & $\overline\ell$ & $0.12$ \\
\cline{2-3} 
          & $\sigma$ & $0.20$ \\
\hline
\multirow{2}{*}{Hybrid}   & $\overline\ell_s$, $\sigma_s$ & $0.12$, $0.20$ \\
\cline{2-3} 
            & $\overline\ell_{\overline r}$, $\sigma_{\overline r}$, $\nu_{\overline r}$ & $0.05$, $0.12$ , $1.5$  \\                
\cline{2-3} 
            & $n_{\mathrm{rect}}$ & $5.0$ \\ 
\cline{2-3} 
            & $\left[w_\mathrm{\min},\, w_\mathrm{\max}\right]$ & $\left[0.02, 0.10\right]$ \\
\cline{2-3} 
            & $\left[h_\mathrm{\min},\, h_\mathrm{\max}\right]$ & $\left[-0.35, 0.35\right]$ \\
\hline
\end{tabular}
\caption{Parameters for three random area generators.}
\label{tab:params}
\end{table}

\begin{figure}[htp]
    \centering
    \includegraphics[width=\linewidth]{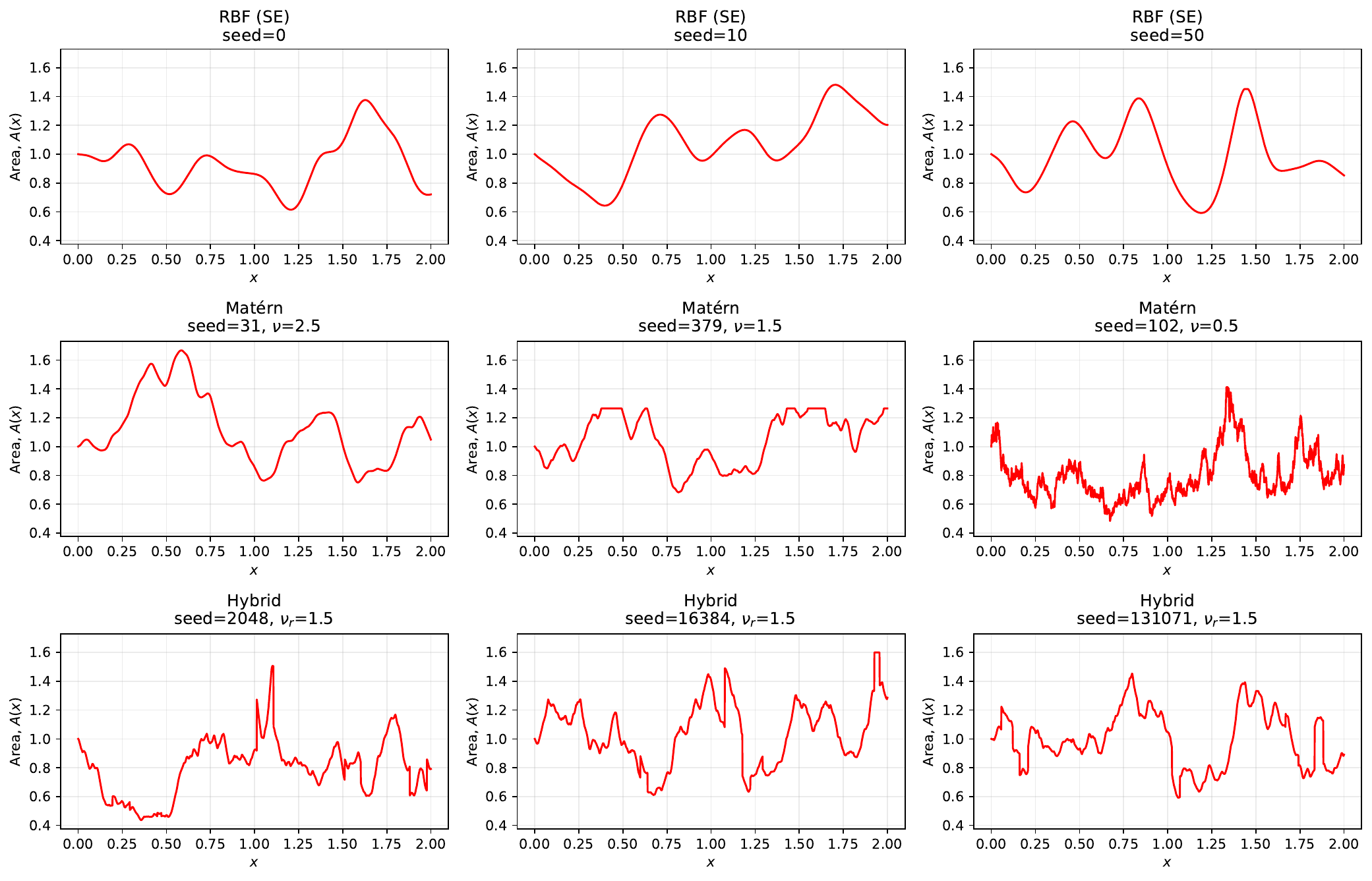}
    \caption{Example area profiles: squared-exponential (top), Mat\'ern (middle) and hybrid model (bottom).}
    \label{fig:area_samples}
\end{figure}

The generation procedure (Algorithm~\ref{alg:area_generation}) computes the Cholesky factor $\mathbf{L}$ of $\mathbb{K}$, draws $\mathbf{z}\sim\mathcal{N}(\mathbf{0},\mathbb{I})$, sets $\mathbf{g}=\mathbf{L}\mathbf{z}$, maps to area via $A_0\exp(\mathbf{g})$, normalises and clips. Example realisations appear in Figure~\ref{fig:area_samples}. Here $\mathbb{I}$ denotes the identity matrix. Since
$\mathbf{L}\mathbf{L}^{\top}=\mathbb{K}$ and
$\mathbf{z}\sim\mathcal{N}(\mathbf{0},\mathbb{I})$, the vector
$\mathbf{g}=\mathbf{L}\mathbf{z}$ has covariance $\mathbb{K}$ as required.

\begin{algorithm}[H]
\begin{algorithmic}[1]
\REQUIRE Spatial grid $\{x_i\}_{i=0}^N$, kernel choice and parameters, baseline  area $A_0$, bounds $[A_{\min}, A_{\max}]$
\ENSURE  Area profile $\mathbf{A} = [A(x_0),\dots,A(x_N)]^\top$
\STATE   Construct $\mathbb{K}$; add jitter $\varepsilon \mathbb{I}$; compute Cholesky factor $\mathbf{L}$ 
\STATE Draw $\mathbf{z} \sim \mathcal{N}(\mathbf{0},\mathbb{I})$; set $\mathbf{g} = \mathbf{L}\mathbf{z}$ 
\STATE Map to area: $\mathbf{A} = A_0 \exp(\mathbf{g})$; normalize: $\mathbf{A} \leftarrow \mathbf{A} / A(0)$; clip
\RETURN $\mathbf{A}$
\end{algorithmic}
\caption{Random area profile generation algorithm}
\label{alg:area_generation}
\end{algorithm}

\section{Numerical implementations}
\label{section3:implementation}

\subsection{Data generation (direct problem simulation)}

\subsubsection{Staggered Leapfrog discretisation of the First-Order Webster system}\label{SG_directpart}

The first-order system \eqref{eq:intro-wave-eq} is solved on $x\in [0,\ell]$ using a staggered-grid leapfrog scheme \cite{schneider2010understanding}. The spatial domain is discretised into $J$ intervals of equal length $\Delta x.$ The flow velocity $Q$ is placed at the nodes 
\begin{equation}
x_j=j\Delta x,\qquad j=0,\cdots,J,
\end{equation}
while the pressure head $H$ is placed at the half-nodes 
\begin{equation}
x_{j+\frac{1}{2}}=\left(j+\frac{1}{2}\right)\Delta x,\qquad j=0,\cdots,J-1.    
\end{equation}
Time is discretised with step $\Delta t$ and the staggered leapfrog updates are

\begin{align} 
Q_j^{\,n+1} &= Q_j^{\,n} - \frac{\Delta t}{\Delta x}\, A_j \bigl(H_{j+\frac{1}{2}}^{\,n+\frac{1}{2}} - H_{j-\frac{1}{2}}^{\,n+\frac{1}{2}}\bigr), \quad \quad\quad\,\, j=1,\cdots,J-1, \label{leap_frogQ}\\
H_{j+\frac{1}{2}}^{\,n+\frac{3}{2}} &= H_{j+\frac{1}{2}}^{\,n+\frac{1}{2}} - \frac{\Delta t}{\Delta x}\, \frac{1}{A_{j+\frac{1}{2}}} \bigl(Q_{j+1}^{\,n+1} - Q_j^{\,n+1}\bigr), \quad j=0,\cdots,J-1. \label{leap_frogH}
\end{align}
Here $A_j=A(x_j)$ and the half-node area is obtained by linear interpolation
\begin{equation}
A_{j+\frac{1}{2}}= \frac{1}{2}\left(A_j + A_{j+1}\right).
\end{equation}
The staggered central differences give a formally second-order approximation in space, and the leapfrog update is second-order in time. Hence, for sufficiently smooth area profiles, the global discretisation error is of order

\begin{equation}
\mathcal{O}\left(\Delta x^2+\Delta t^2\right).\label{eq:order}     
\end{equation}
The boundary data used in the cross-transfer experiments of Section 4.3 are generated by this scheme, so the simulated measurements are themselves accurate only to the order given in \eqref{eq:order}, and that limit carries through to the reconstruction errors reported there. A refinement study quantifying its size is given in Subsection \ref{sec:forward-convergence}.

In the computations reported here, we take \(\Delta t=\Delta x\), so that the Courant number is one and the leapfrog scheme is non-dissipative in the standard sense \cite{strikwerda2004finite}. The waveguide is excited at the inlet through the left boundary value of the flow variable. We use the discrete unit impulse
\begin{equation}
Q_0^{n+1}=\begin{cases}\label{discrete_delta}
    \frac{1}{\Delta t}\,: \quad\, n=0,\\
    0 \quad: \quad n\geq 1,
\end{cases}
\end{equation}
which satisfies $\sum_{n}u_0^{n+1}\Delta t=1$. This provides a normalized discrete approximation of a Dirac delta source at the inlet. We also tested a temporally spread pulse,
\begin{equation}
Q_0^{n+1} =
\begin{cases}\label{discrete_spread_delta}
\frac{1}{K_\text{imp}\Delta t}\,: \quad\,\, n=0,\cdots,K_\text{imp}-1,\\
0 \quad\quad\,\,\,\,\,: \quad\,\,\, n\geq K_\text{imp},
\end{cases}
\end{equation}
where \(K_{\mathrm{imp}}\) denotes the number of time steps over which the impulse is distributed. This did not produce noticeable changes in the computed responses or reconstructions. Therefore, unless otherwise stated, all SG results use the one-step discrete impulse \eqref{discrete_delta}.

To reduce the non-physical reflections from the truncated computational domain $[0,\ell],$ we impose a first-order Engquist-Majda absorbing boundary condition at the right endpoint $x=\ell.$ For an outgoing wave of the form $u(x,t)\approx F(t-x/c_0)$, the boundary satisfies the one-way equation 
\begin{equation}
u_t + c_0u_x=0 \text{ at } x=\ell. \label{absorbing_boundary}
\end{equation}
Let $\lambda=c_0\Delta t/\Delta x$ denote the Courant number and let $x_J=\ell$ be the last grid point. Discretising \eqref{absorbing_boundary} at the right boundary gives the update
\begin{equation}
Q_J^{n+1} \approx Q^{n}_{J-1}+\frac{1-\lambda}{1+\lambda}\left(Q^n_J -Q_{J-1}^{n+1}\right).\label{right_boundary_discrete}
\end{equation}
This is the update used after the interior values have been advanced. Since the computations for this section take $\Delta t = \Delta x$, we have $\lambda=1$ and \eqref{right_boundary_discrete} reduces to $Q_J^{n+1}=Q_{J-1}^n.$

\subsubsection{Central Leapfrog discretisation of the Second-Order Webster equation}\label{KLO_directpart}
The KLO forward data are generated by solving the second-order Webster equation
\begin{equation}
\frac{\partial^2 }{\partial t^2}H(x,t) = \frac{c^2_0}{A(x)}\frac{\partial}{\partial x} \left(A(x)\frac{\partial}{\partial x}H(x,t)\right), \quad x\in [0,\ell],\, t>0,\label{KLO_direct}
\end{equation}
with zero initial conditions $H(x,0)=0,$ and $H_t(x,0)=0$. At the left boundary $x=0,$ the waveguide is excited through the Neumann condition, \(H_x(0,t)=f(t).\) At the right end, $x=\ell$, we impose an absorbing condition to reduce artificial reflections from the computational boundary. The spatial interval is discretized using the uniform grid
\begin{equation}
x_i=i\Delta x,\quad i =0,1,\cdots, n_x-1,\quad \Delta x = \frac{\ell}{n_x-1},\label{KLO_spatial}   
\end{equation}
and time interval by 
\begin{equation}
t_n=n\Delta t,\quad n =0,1,\cdots, n_t-1.\label{KLO_time}
\end{equation}
The area is sampled as $A_i=A(x_i),$ and the midpoint values are approximated by
\begin{equation*}
A_{i+\frac{1}{2}}=\frac{1}{2}\left(A_{i+1}+A_i\right).    
\end{equation*}
To approximate the spatial operator in conservative form, we first define the half-grid flux
\begin{equation}
F_{i+\frac{1}{2}}^n = A_{i+\frac{1}{2}}\frac{H^{n}_{i+1}-H^n_i}{\Delta x}.   
\end{equation}
The divergence term in \eqref{KLO_direct} is then approximated by
\begin{equation}
\frac{F_{i+\frac{1}{2}}^n - F_{i-\frac{1}{2}}^n}{\Delta x}.    
\end{equation}
Thus, the method first approximates \(H_x\) on half-grid points, forms the flux \(AH_x\), and differentiates this flux back to the nodal grid. For constant \(A\), this reduces to the standard second-difference approximation of \(H_{xx}\).
The time derivative is discretised by the central three-level leapfrog formula
\begin{equation}
H_{tt}(x_i,t_n) = \frac{\left(H^{n+1}_i-2H^n_i+H^{n-1}_i\right)}{\Delta t^2}.\label{KLO_time_discrete}
\end{equation}
For interior nodes \(i=1,\ldots,n_x-2\), the update becomes
\begin{equation}
H^{n+1}_i = 2H^n_i-H^{n-1}_i + \frac{\Delta t^2}{\rho_i \Delta x} \left(F_{i+\frac{1}{2}}^n-F_{i-\frac{1}{2}}^n\right), \quad \rho_i =\frac{A_i}{c_0^2}.
\end{equation}
The scheme is explicit, and the time step is chosen according to the CFL condition $\lambda:=\frac{c_0\Delta t}{\Delta x} < 1,$ where $\lambda$ denotes the Courant number. In the numerical experiments, we use $\Delta t = \eta \frac{\Delta x}{c_0},\quad 0<\eta < 1,$ with final computations performed using \(\lambda=0.4\).
 
At the right boundary, we use an Engquist--Majda type absorbing update,
\begin{equation}
H^{n+1}_{n_x-1} = H^{n}_{n_x-2} + \frac{\lambda -1}{\lambda +1}\left(H^{n+1}_{n_x-2} - H^{n}_{n_x-1}\right).\label{right_boundary}
\end{equation}
In the case with $\lambda=1,$  \eqref{right_boundary}  reduces to 
\begin{equation}
H^{n+1}_{n_x-1} = H^{n}_{n_x-2},
\end{equation}
corresponding to a one-grid-step transport of the outgoing wave. In practice, the modified absorbing update was found to give better right-end behaviour than the simpler first-order outflow update, particularly for smaller Courant numbers. 

At the left boundary $x=0$, the Neumann excitation is imposed through the boundary flux 
\begin{equation*}
F^{n}_{-\frac{1}{2}}=A_0 f(t_n),
\end{equation*}
where $f(t_n)$ is taken as a discrete approximation of the Dirac delta input. The input \(f(t_n)\) is taken as a discrete approximation of a Dirac impulse,
\begin{equation}
f{(t_n)}=\begin{cases}\label{discrete_delta_for_KLO}
\frac{1}{\Delta t}\,: \quad\, n=0,\\
    0 \quad: \quad n\geq 1.
\end{cases}
\end{equation}
The interior-side flux at the first half-grid point is
\begin{equation*}
F^{n}_{\frac{1}{2}}=A_{\frac{1}{2}}\frac{H^n_1-H^n_0}{\Delta x},
\end{equation*}
The first grid point is then updated by a conservative half-cell balance over the interval $[0,\Delta x/2]$, giving
\begin{equation}
H^{n+1}_0=2H^n_0 -H^{n-1}_0+ \frac{\Delta t^2}{\rho_0}\frac{F^n_{\frac{1}{2}}-F^n_{-\frac{1}{2}}}{\Delta x/2},\quad \rho_0 =\frac{A_0}{c^2_0}.\label{boundary_flux}
\end{equation}
This incorporates the prescribed Neumann input without introducing a ghost point and is consistent with the conservative form of the Webster equation. The simulated measurement is the left boundary trace
\begin{equation}
H^{n}_0\approx H(0,t_n)\quad n =0,\cdots,n_t-1.\label{KLO_impulse_response}
\end{equation}
This trace is used as the boundary impulse response for the KLO inverse reconstruction described in Subsection~\ref{KLO_inverse}.

\subsection{Algorithm Inversion}

\subsubsection{Sondhi--Gopinath inversion algorithm} \label{SG_inverse}
Following \cite{sondhi71_deter_vocal_shape_from_impul}, the relationship between area and the impulse response is attained as a Fredholm integral equation of the second kind \eqref{eq:fredholm}  where $h_\mathrm{SG}\left(|t-s|\right)$ denotes the continuous or theoretical impulse response kernel. The area profile is the recovered from the end point value of $f$ via \eqref{eq:area-from-f-squared}.

To solve \eqref{eq:fredholm} numerically, we discretise the integral equation on the temporal grid $t_n=n\Delta t.$  For a reconstruction parameter $a=m\Delta t$, the system is 
\begin{equation}
\left(\mathbb{I}+\mathbf{K}(a)\right)\mathbf{f}=\mathbf{1}_M,
\end{equation}
where the $(2m+1)\times(2m+1)$ matrix $\mathbf{K}(a)$ is Toeplitz whose kernel depends only on $|i-j|$. We add a small regularisation parameter $\varphi >0$ to stabilise the inversion and solve \(\mathbf{T}(a)\mathbf{f}=\mathbf{1}_M\) with
the Toeplitz matrix $\mathbf{T}(a)$ defined by its first column as
\begin{equation}
c_k=\frac{\Delta t}{2}{\mathbf{\tilde h}}(k\Delta t)+(1+\varphi)\delta_{k0},\quad k = 0,1,\cdots,M-1,
\end{equation}
using the Levinson algorithm via \texttt{scipy.linalg.solve\_toeplitz} ($O(M^2)$). Here $\delta_{k0}$ denotes the Kronecker delta. 

To improve the quadrature accuracy, we replace the trapezoidal rule used in the discretisation by fourth-order Gregory weights \cite{fornberg2019improved, ralston2001first}. For $M\geq 5$, the Gregory weights are
\begin{equation}
w_j = \Delta t
\begin{cases}
\frac{5}{12}, \,\,\, j= 0 \text{ or } J=M-1,\\
\frac{13}{12}, \,\,\, j= 1 \text{ or }  J=M-2,\\
1,\quad \mathrm{otherwise}.\\
\end{cases}
\end{equation}
The correction is a rank-$p$ ($p\leq4$) perturbation handled via the Woodbury identity \cite{golub2013matrix}. The solution vector $\mathbf{f}$ approximates the values of $f(a,\tau)$ on the grid 
\begin{equation}
\tau_j=-a+j\Delta t, \qquad j =0,\cdots,M-1.
\end{equation}
Hence, using \eqref{eq:area-from-f-squared}, the reconstructed area at $x=a$ is obtained from the endpoint value
\begin{equation}
A_\mathrm{rec}(a)=f^2_{M-1}.  
\end{equation}
\subsubsection{Korpela--Lassas--Oksanen inversion algorithm}\label{KLO_inverse}
The KLO inverse algorithm reconstructs the cross-sectional area \(A(x)\) from the measured boundary impulse response  $\mathbf{\hat h}:=\mathbf{\hat h(0,t)}$. The method follows the discrete regularised formulation of \cite{korpela16_regul_strat_inver_probl_1}, in which the measured data are used to construct a connecting operator $\mathbf{K}$. Solving a family of regularised linear systems over increasing time windows then yields a scalar function whose derivative determines the area.

Let the time interval be discretised uniformly as
\begin{equation}
t_i=i\Delta t, \quad \text{ for } i=0,\cdots, n_t-1, 
\end{equation}
and let the measured impulse response be represented by the vector $\mathbf{\hat h}\in  \mathbb{R}^{n_t}$. In the implementation, the vector used in the convolution operator is the scaled discrete impulse response
\begin{equation}
\hat{h}_i=\Delta t H(0,t_i)\label{scaled_discrete_impulse}
\end{equation}
The elementary operators are as follows. 
\begin{enumerate}
    \item \textbf{Time reversal matrix} $\mathbf{R}$. The operator $\mathbf{R}\in\mathbb{R}^{n_t\times n_t}$ represents time reversal on the discrete time grid. For a vector $\mathbf{v}=(v_0,v_1,\cdots,v_{n_t-1})^\top$, its action is defined by
    \begin{equation}
    \mathbf{Rv}_i=v_{n_t-1-i}, \qquad i=0,\cdots,n_t-1.   
    \end{equation}
    Consequently, $\mathbf{R}$ is anti-diagonal identity matrix. As a symmetric orthogonal permutation matrix, it inherently satisfies
    $\mathbf{R}^\top=\mathbf{R}^{-1}=\mathbf{R}$ and $\mathbf{R}^2=\mathbb{I}$.
    \item \textbf{Integration matrix} $\mathbf{J}.$ The operator $\mathbf{J}$ approximates the integral  operator 
    \begin{equation}
    Jf(t)=\frac{1}{2}\int_t^{2T_0-t} f(s)ds,\qquad 0\leq t\leq T_0,\label{integral_operator_J}
    \end{equation}
    where $T_0=t_{n-1}/2.$ On a uniform grid, this is represented by the matrix
    \begin{equation}\label{operator_J}
     \mathbf{J}_{i,\overline\tau} = 
    \begin{cases}
    \frac{\Delta t}{2}, \quad i\leq \overline{\tau}\leq N-1 \text{ and } N-i>i,\\
    0, \quad \quad \mathrm{otherwise},
    \end{cases}    
    \end{equation}
    for $i,\overline{\tau}=0,\cdots,N.$
    \item \textbf{Convolution matrix} $\mathbf{\Lambda}.$ 
   For a linear time-invariant system, the mapping from the source signal $f(t)$ to the boundary input ${H(0,t)}$ is given by convolution with the impulse response $\hat{h}(t).$  On the uniform grid $t_i=i\Delta t$, this Neumann-to-Dirichlet map is represented by a lower-triangular Toeplitz convolution matrix
  $\mathbf{\Lambda}\in\mathbb{R}^{n_t\times n_t}$
   defined by 
   \begin{equation}\label{convolution_matrix}
   \mathbf{\Lambda}_{i,j}=
   \begin{cases}
   \hat h_{i-j}, \text{ for } i\geq j, \\
   0,\quad\,\,\, \text{ for } i<j,
   \end{cases}
   \qquad i,j =0,\cdots,n_t-1.
   \end{equation}
   Here, the faction $\Delta t$ from the discrete convolution has already been absorbed into the impulse response vector in \eqref{scaled_discrete_impulse}. Equivalently, one could unscale \eqref{scaled_discrete_impulse} and instead write the matrix entries as
   \begin{equation}
   \mathbf{\Lambda}_{i,j}=
   \begin{cases}
   \hat h_{i-j}\Delta t, \text{ for } i\geq j, \\
   0,\quad\,\,\, \text{ for } i<j,
   \end{cases}
   \qquad i,j =0,\cdots,n_t-1.
   \end{equation}
    
   \item \textbf{The} $\mathbf{K}$
   \textbf{operator.} The KLO connecting operator is then assembled as¨
   \begin{equation}\label{K_matrix}
    \mathbf{K}=\mathbf{R\Lambda R J}-\mathbf{J\Lambda}.
    \end{equation} 
    To compute the regularised solution for a fixed time window, let \(T_0=t_{n_t-1}/{2}.\) For each parameter $r\in [0,T_0],$ define the time window 
    \begin{equation}
    \mathcal{I}_r:=[T_0-r,T_0],
    \end{equation}
    and let  $\mathcal{J}_r:=\left\{i:t_i\in \mathcal{I}_r\right\}$ be the corresponding index set.  The matrix \(K\) is restricted to this window by extracting
    \begin{equation}
    \mathbf{K}_{rr}=\mathbf{K}[\mathcal{J}_r,\mathcal{J}_r],   
    \end{equation}
    that is, the restriction of $\mathbf{K}$ to the rows and columns corresponding to the time window $\mathcal{I}_r$. 
\end{enumerate}

The right-hand side is the discretisation of
\begin{equation}
B_1(t)=\begin{cases}\label{operator_b1}
    T_0-t, \quad 0\leq   t\leq T_0,\\
    0, \quad\quad\quad   t> T_0,
\end{cases}
\end{equation}
We represent \eqref{operator_b1} on the discrete time grid by a vector $\mathbf{b}_1\in \mathbb{R}^{n_t}$ with components 
\begin{equation}
\mathbf{b}_{1,i} = B_1(t_i), \qquad i = 0,\cdots, n_t-1. \label{b1_component}
\end{equation}
For a given regularization parameter $\alpha >0$, we solve the linear system restricted to the window $\mathcal{I}_r$
\begin{equation}
\left(\mathbf{K}_{rr}+ \alpha \mathbb{I}\right)\mathbf{f}_{r,\alpha}=\mathbf{b}_{1|\mathcal{J}_r}, \label{KLO_Linear_System}
\end{equation}
where $\mathbf{b}_{1|\mathcal{J}_r}$ denotes the restriction of $\mathbf{b}_{1}$ to the index set $\mathcal{J}_r.$ The solution vector $\mathbf{f}_{r,\alpha}$ is then embedded into a full-length vector $\mathbf{f}_{\mathrm{full}}\in \mathbb{R}^{n_t}$ by assigning its entries on $\mathcal{J}_r$ to $\mathbf{f}_{r,\alpha}$ and setting all remaining entries to  zero. The scalar quantity 
\begin{equation}
s_\alpha(r)=\langle \mathbf{f}_{\mathrm{full}}, \mathbf{b}_{1} \rangle\Delta t,
\end{equation}
is the computed. Repeating this procedure for over a discrete set of $r=r_j$ yields the discrete functions $s_\alpha (r_j).$  Since the restricted systems are ill-conditioned, particularly for small time windows, Tikhonov regularisation is used. Following \cite{korpela16_regul_strat_inver_probl_1}, we set
\begin{equation}
\alpha = \beta\,\varepsilon^{\frac{4}{9}},\qquad\varepsilon = 10^{-4},\qquad \beta = 2\times10^{-5},\label{KLO_regularisation}
\end{equation}
where \(\varepsilon\) is used as a small stabilising parameter rather than as a physical noise level. The area-related quantity is obtained by differentiating \(s_\alpha(r)\). On a uniformly spaced grid $r_j$ with spacing $\Delta r$, we approximate the derivative by 
\begin{equation}
k_0:=\frac{s_1-s_0}{\Delta r}, \quad k_j:=\frac{s_{j+1}-s_j}{\Delta r}\,\,(1\leq j\leq m-2), \quad k_{m-1}:=\frac{s_{m-1}-s_{m-2}}{\Delta r}.\label{KLO_derivative}
\end{equation}
The reconstructed area is then computed from
\begin{equation}
x_j=c_or_j, \qquad A_j=c_0k_j.\label{interpolated_values}   
\end{equation}
The values $A_j$ are clipped to enforce physical bounds, $A_\mathrm{min}\leq A(x)\leq A_\mathrm{max}$. These pairs $(x_j,A_j)$ are then interpolated onto a fine spacial grid $\tilde x\in [0, \ell].$ The reconstruction procedure contains two numerically sensitive steps that may amplify noise and discretization errors: the discrete differentiation \eqref{KLO_derivative} and the interpolation of the discrete pairs \((x_j,A_j)\) onto a finer spatial grid. To reduce these effects, we apply Gaussian smoothing. The smoothing step is implemented as a discrete convolution with a normalised Gaussian kernel \cite{lindeberg2024discrete}. For a uniformly sampled sequence \(y_j=y(\xi_j)\), its smoothed version is defined by
\begin{equation}
\tilde y_j=\sum_{m=-M}^{M}w_m y_{j+m},
\end{equation}
where
\begin{equation}
w_m=\frac{1}{W}\exp\left(-\frac{m^2}{2\sigma_{\mathrm{pts}}^2}\right),
\qquad
W=\sum_{m=-M}^{M}\exp\left(-\frac{m^2}{2\sigma_{\mathrm{pts}}^2}\right).
\end{equation}
The normalisation ensures that \(\sum_m w_m=1\), so that constant signals are preserved. The kernel half-width is chosen as
\begin{equation}
M=\max\left(1,\lfloor 4\sigma_{\mathrm{pts}}\rfloor\right),
\end{equation}
which keeps the kernel compact while retaining most of the Gaussian mass \cite{lindeberg2002scale}. Near the endpoints, the convolution is evaluated using symmetric reflection padding \cite{strang1996wavelets}. This avoids introducing artificial constant extensions at the boundaries and reduces smoothing-induced endpoint artefacts, which is important because the reconstructed area near \(x=0\) and \(x=\ell\) is part of the error evaluation.

After computing the discrete derivative \(k_j\), we apply this Gaussian smoothing to obtain \(\tilde k_j\). This step is used because numerical differentiation can amplify high-frequency oscillations. In the reported experiments, the smoothing parameter is set to \(\sigma_{k,r}=5\), which was found to provide a suitable balance between noise suppression and preservation of the main features of \(k(r)\). The smoothed derivative is then used in the reconstruction formula for \(A(x)\). After interpolation onto the final spatial grid, no additional smoothing is applied to the reconstructed area profile, that is, \(\sigma_{A,x}=0\).

\subsection{Cross algorithm validation}
In this subsection, we describe the cross-algorithm validation used to compare the two methods. Specifically, each algorithm is used to generate data that are then supplied to the other algorithm, allowing bidirectional testing of their reconstruction behaviour.

\subsubsection{SG scheme driven by KLO boundary data}
In its standard formulation, the SG reconstruction depends on the convolution between the inlet flow  $Q(0,t)$ and the inlet  pressure head $H(0,t)$ discussed in \eqref{convolution_of_SG}.  In our coupled framework, the SG algorithm bypasses its own internal forward model. Instead, the required impulse response kernel is extracted from the boundary trace generated by the KLO forward simulation. 

We connect the SG and KLO data by applying a Neumann impulse at the inlet. This impulse induces a step-type inlet flow, scaled by the inlet area, so that the corresponding inlet head is the convolution of the SG impulse response with this step input
\begin{equation}
H(0,t)=-A(0)\int_0^t h_{\mathrm{SG}}(t-s) ds.\label{KLO_SG_Relations}
\end{equation}
The KLO boundary trace determines the SG impulse-response kernel through
\begin{equation}
h_{\mathrm{SG},\mathrm{KLO}}(t) = -\frac{1}{A(0)}\frac{d}{dt} H(0,t),\quad t\geq 0.\label{hSG_from_KLO}
\end{equation}
The resulting discrete kernel is then inserted directly into the SG inversion procedure.

In the noisy experiments, the perturbation is introduced only in the KLO-generated inlet boundary trace, as described in \eqref{sigma_KLO} and \eqref{Noisyimpulse_KLO}. The same noisy trace \eqref{Noisyimpulse_KLO} is then used consistently in both reconstructions. On the KLO side, it defines the noisy kernel used to construct \(\mathbf{\Lambda}\) and hence \(\mathbf{K}\). On the SG side, a smoothed version of \eqref{Noisyimpulse_KLO} is differentiated to obtain the corresponding discrete SG impulse response.

\subsubsection{KLO scheme driven by SG boundary data}\label{sec:klo-from-sg}
In this case, the KLO Neumann-to-Dirichlet map is not generated by a separate KLO forward solver, but is approximated from the SG inlet impulse response. The connection follows from the Webster relation \(\partial_t Q(x,t)=-A(x)\partial_x H(x,t),\) which gives, at the inlet,
\begin{equation}
Q(0,t)=-A(0)\int_0^t f(s)\,ds,\quad f(t)=H_x(0,t).
\end{equation}
Thus, the conversion from the SG flow-driven data to the KLO Neumann-driven map requires the time-integrated SG impulse-response kernel. If \(h_{\mathrm{SG}}\) denotes the normalised SG flow-to-head impulse response, then
\begin{equation}
H(0,t) = -A(0)\int_0^t\left(\int_0^{t-s} h_{\mathrm{SG}}(\xi)\,d\xi
\right)f(s)\,ds.
\end{equation}
In the discrete setting, the inner time integration is implemented by cumulative sums of the SG impulse-response samples.

On a uniform temporal grid $t_n=n\Delta t$, let $h_{\mathrm{SG,n}}$ denote the normalised SG kernel samples and define
\begin{equation}
c_n = \sum_{m=0}^nh_{\mathrm{SG,m}}.\label{eq:cumsum}
\end{equation}
The discrete Neumann-to-Dirichlet map used in the KLO reconstruction is then approximated by the lower-triangular Toeplitz matrix
\begin{equation}
(\mathbf{\Lambda}_{\mathrm{SG}})_{nk} =
\begin{cases}
-A_0\,\Delta t^2\,c_{n-k}, & n\ge k,\\[4pt]
0, & n<k,
\end{cases}
\qquad n,k=0,\cdots,N-1.
\end{equation}
Once \(\mathbf{\Lambda}_{\mathrm{SG}}\) is assembled, the remaining area reconstruction proceeds as described in Subsection~\ref{KLO_inverse}. 

The two transfer directions are not equivalent in their sensitivity to perturbation. Passing KLO data to the SG algorithm requires the
boundary trace to be differentiated, \eqref{hSG_from_KLO}, whereas passing SG data to the KLO algorithm requires it to be integrated, \eqref{eq:cumsum}. Differentiation is unbounded: on a grid of step $\Delta t$ a perturbation of size $\varepsilon$ produces an error of order $\varepsilon/\Delta t$ in the
derivative, so high-frequency components of the perturbation are amplified. Cumulative summation is bounded and attenuates them, and uncorrelated errors
partially cancel. The implementation reflects this: in the KLO-to-SG direction the trace must be smoothed before the derivative is taken, whereas the SG-to-KLO direction requires no such regularisation of the data. This is why the ensemble study of Subsection~\ref{sec:statistical-validation} is carried out in the SG-to-KLO direction, while the reverse direction is illustrated by representative examples only.

\subsection{Noise case study for SG}\label{Noise_4_SG}

The noise-free SG impulse response is obtained using the forward model described in Subsection \ref{SG_directpart}. The forward simulation produces the inlet trace on the forward time grid. This trace is then resampled onto the inverse time grid. Before adding noise, the direct contribution of the imposed boundary impulse is removed. The processed noise-free impulse response used in the inverse reconstruction is thereby defined by \(\mathbf{\tilde h}:=\left(\tilde h_0, \tilde h_1, \cdots, \tilde h_{N-1}\right)^\top.\)

This is the signal to which noise is added in the numerical experiment. The signal was then corrupted by additive Gaussian noise at the relative levels
\begin{equation}
\delta\in\{0.0,\,0.01,\,0.05,\,0.10\}.\label{noise_levels}
\end{equation}
Let \(\tilde\eta = \left(\tilde\eta_0,\cdots,\tilde\eta_{N-1}\right)^\top\) be a vector of independent standard Gaussian random variables, $\tilde\eta_n\sim\mathcal{N}(0,1).$ 
A fixed noise direction $\eta^0$ is then obtained by
\begin{equation}
\eta^0:=\frac{\tilde \eta}{\|\tilde \eta\|_{2,\Delta t}}, \quad \|\eta^0\|_{2,\Delta t} =1,
\end{equation}
with $\|\cdot\|_{2,\Delta t}$ as the weighted discrete $L^2$-norm. For a prescribed relative noise level $\delta$, the perturbation is scaled by the energy of the noise-free signal,
\begin{equation}
\eta^\delta := \delta\|\mathbf{\tilde h}\|_{2,\Delta t}\, \eta^0,\qquad \|\eta^\delta\|_{2,\Delta t} = \delta\|\mathbf{\tilde h}\|_{2,\Delta t}.\label{SG_noiseimpulse}
\end{equation}
The noisy impulse response is then defined by \(\mathbf {\tilde h}^\delta = \mathbf{\tilde h} + \eta^\delta. \)

\subsection{Noise case study for KLO}\label{Noise_4_KLO}
We consider the Neumann-to-Dirichlet map \(\Lambda: f(t) \rightarrow H(0,t)\)
associated with the second order Webster system in \eqref{eq:intro-wave-eq}. In practice, only a perturbed version $\Lambda^\varepsilon$ is available, satisfying
\begin{equation}
\|\Lambda^\varepsilon -\Lambda\| \leq \varepsilon,
\end{equation}
where $\|\cdot\|$ denotes the operator norm. In regularized reconstruction procedure, the regularization parameter is chosen according to \cite{korpela16_regul_strat_inver_probl_1} and \eqref{KLO_regularisation}. Here $\beta >0$ is selected empirically depending on the noise regime. The Neumann-to-Dirichlet map is represented through its discrete impulse response kernel \(\mathbf {\hat h} := \left(\hat h_0, \hat h_1,\cdots ,\hat h_{n_t-1} \right)^{\top},\) obtained from the boundary impulse experiment described in Subsubsection \ref{KLO_directpart}. The corresponding discrete Neumann-to-Dirichlet is a Toeplitz matrix whose size is measured by an estimate of the induced operator norm $\|\mathbf\Lambda\|_{2}$, computed by power iteration using the Toeplitz convolution structure. For the noisy experiments, we consider \eqref{noise_levels}. For $\delta>0$, the effective perturbation size is defined by
\begin{equation}
\varepsilon_\mathrm{eff}(\delta) = \delta \|\mathbf{\Lambda}\|_2.\label{sigma_KLO}
\end{equation}
A single random noise realization $n_h$ is drawn once using a fixed seed. The noise vector is then normalized so that the Toeplitz operator generated by it has unit estimated operator norm. Denoting this fixed normalized noise direction by $\mathbf{\hat  n}_h$, the perturbed impulse-response kernel is given by
\begin{equation}
\mathbf{\hat h}^\delta = \mathbf{\hat h} + \varepsilon_\mathrm{eff}(\delta)\mathbf{\hat{n}}_h.\label{Noisyimpulse_KLO}
\end{equation}
The regularization parameter is chosen as $\alpha=\beta \varepsilon_\mathrm{eff}(\delta)^{4/9}$, with the prefactor $\beta$ selected according to the noise level
\begin{equation}
\beta =
\begin{cases}
2\times 10^{-5},\quad \delta=0,\\
5\times 10^{-3},\quad 0< \delta\leq 0.01,\\
1\times 10^{-2},\quad \delta>0.01.
\end{cases}
\end{equation}
Gaussian smoothing width scales as 
\begin{equation}
\sigma_{k,r} =\sigma_{\mathrm{ref}}\frac{N_{\mathrm{grid}}}{N_{\mathrm{ref}}}\bigl(1+40\delta\bigr),
\end{equation}
where \(N_{\mathrm{grid}}\) is the number of points in the current \(r\)-grid, with \(N_{\mathrm{ref}}=1500\), and \(\sigma_{\mathrm{ref}}=6\). 

\subsection{Statistical comparison design}

We perform a paired Monte Carlo study with $n=1000$ independent area realisations. For each realisation, both methods are applied to the same noisy data set at each noise level $\delta$, yielding paired error observations $(e_i^{\mathrm{SG}},e_i^{\mathrm{KLO}})$. To avoid inverse crime \cite{kaipio2005statistical}, the forward problem is solved on a fine spatial grid and data are resampled to a coarser measurement grid before inversion.

\subsubsection{Error metrics}\label{error_measure}

To compare the accuracy of the SG and KLO reconstructions in a consistent manner, all errors are evaluated on a common spatial grid. Let $A_\mathrm{true}(x)$ denote the true area profile and $A_\mathrm{rec}(x)$ a reconstructed profile. Let \(e(x)=A_\mathrm{true}(x)-A_\mathrm{rec}(x),\) we report absolute and relative $L^2$- and $H^1$-errors
\begin{equation}
\|e\|^\mathrm{rel}_{L^2} = \frac{\|e\|_{L^2}}{\|A_{\mathrm{true}}\|_{L^2}},\quad \|e\|^\mathrm{rel}_{H^1} = \frac{ \left(\|e\|^2_{L^2}+ \|e'\|^2_{L^2}\right)^{1/2}}{\|A_{\mathrm{true}}\|_{H^1}}.
\end{equation}
For each realisation model, noise level, and reconstruction method, we therefore record the four quantities
\begin{equation}
\|e\|^{\mathrm{abs}}_{L^2}, \qquad \|e\|^\mathrm{rel}_{L^2}, \qquad\|e\|^{\mathrm{abs}}_{H^1},\qquad \|e\|^\mathrm{rel}_{H^1}.
\end{equation}

\subsubsection{Hypothesis tests}
For a given error, the paired difference for each realisation $i$ is
\begin{equation}
d^{e}_i = e_i^{\text{KLO}}-e_i^{\text{SG}}, \qquad i=1,\cdots,n.\label{difference_between_SG_and_KLO}
\end{equation}
We test \(H_0:\mathbb{E}[d^e]=0\) using both both a paired $t$-test and the Wilcoxon signed-rank test \cite{wilcoxon1945individual}. Effect sizes include Cohen's $d_{\mathrm{Cohen}}=\bar{d}/s_d$ \cite{cohen1988statistical} and the rank-biserial correlation $r_{\mathrm{rb}}=(W^+-W^-)/(W^++W^-)$.

\section{Results}
\label{section4:results}

\subsection{SG algorithm}
In this subsection, we evaluate the performance of the SG algorithm in two distinct stages. We begin by establishing the method's baseline accuracy in a noise-free setting, quantifying the reconstructions across smooth, rough, and hybrid area profiles using relative and absolute $L^2$- and $H^1$-errors. With this baseline established, we then introduce additive measurement noise to test the algorithm's stability, exploring how varying noise levels and boundary treatments impact the final geometric reconstructions.

\subsubsection{Noiseless case}

\begin{figure}[htbp]
    \centering    \includegraphics[width=\textwidth]{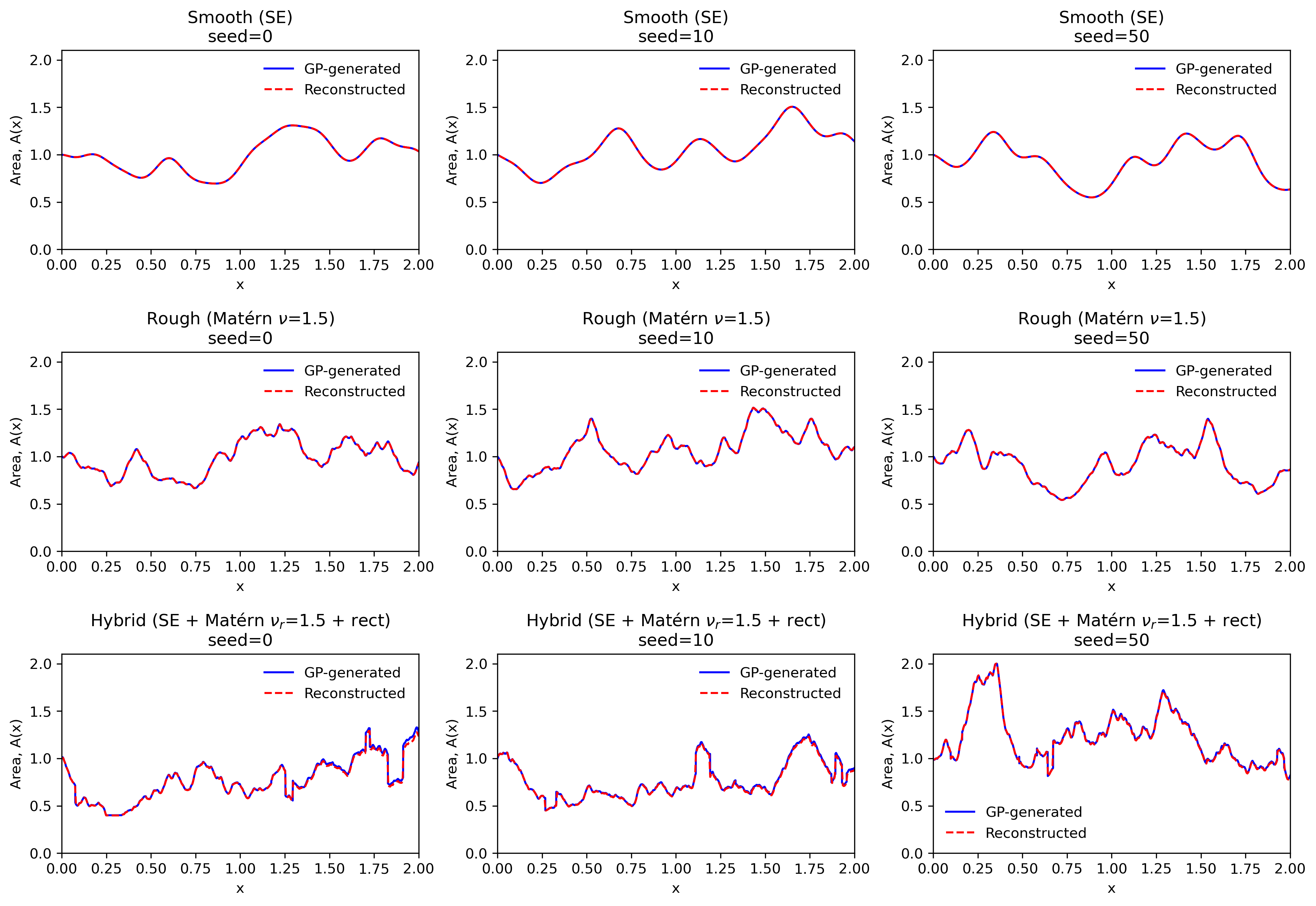}
    \caption{\textbf{Noise-free SG reconstructions.} True profiles (solid blue) vs. reconstructions (dashed red) smooth SE (top), rough Mat\'ern with $\nu=1.5$ (middle), and hybrid (bottom) profiles at seeds $0,10,50$.}
    \label{fig:sg_noiseless_examples}
\end{figure}

Figure~\ref{fig:sg_noiseless_examples} shows noise-free SG reconstructions for the three profile classes (seeds $0,10,50$). The method achieves high accuracy for smooth SE profiles ($L^2_{\mathrm{rel}}\approx5\times10^{-4}$, $H^1_{\mathrm{rel}}\lesssim3\times10^{-3}$) and moderate accuracy for rough Mat\'ern profiles. Hybrid profiles are most challenging, with $L^2_{\mathrm{rel}} \approx 1.8 \times 10^{-2}$ and $H^1_{\mathrm{rel}} \approx 5.5 \times 10^{-2}$, since the $H^1$ metric penalises sharp transitions.  Selected errors are in Table~\ref{tab:sg_noiseless_errors}.

\begin{table}[htbp]
\centering
\begin{tabular}{llcccc}
\toprule
Seed & Case & $H^1_{\mathrm{abs}}$ & $H^1_{\mathrm{rel}}$ & $L^2_{\mathrm{abs}}$ & $L^2_{\mathrm{rel}}$ \\
\midrule
\multirow{3}{*}{0}
  & Smooth & $3.7136\times 10^{-3}$ & $1.6053\times 10^{-3}$ & $2.6894\times 10^{-4}$ & $1.8940\times 10^{-4}$ \\
  & Rough  & $2.9631\times 10^{-2}$ & $7.8039\times 10^{-3}$ & $6.0047\times 10^{-4}$ & $4.2552\times 10^{-4}$ \\
  & Hybrid & $2.8521\times 10^{-1}$ & $3.5189\times 10^{-2}$ & $2.0355\times 10^{-2}$ & $1.7959\times 10^{-2}$ \\
\midrule
\multirow{3}{*}{10}
  & Smooth & $6.4846\times 10^{-3}$ & $2.3587\times 10^{-3}$ & $6.4116\times 10^{-4}$ & $4.2150\times 10^{-4}$ \\
  & Rough  & $3.6418\times 10^{-2}$ & $8.2715\times 10^{-3}$ & $1.1595\times 10^{-3}$ & $7.5329\times 10^{-4}$ \\
  & Hybrid & $3.7151\times 10^{-1}$ & $5.5006\times 10^{-2}$ & $1.3851\times 10^{-2}$ & $1.2438\times 10^{-2}$ \\
\bottomrule
\end{tabular}
\caption{\textbf{Error quantification for noise-free SG reconstructions (selected seeds).}}
\label{tab:sg_noiseless_errors}
\end{table}

\subsubsection{Performance under additive noise}

\begin{figure}[htbp]
\centering
\includegraphics[width=\textwidth]{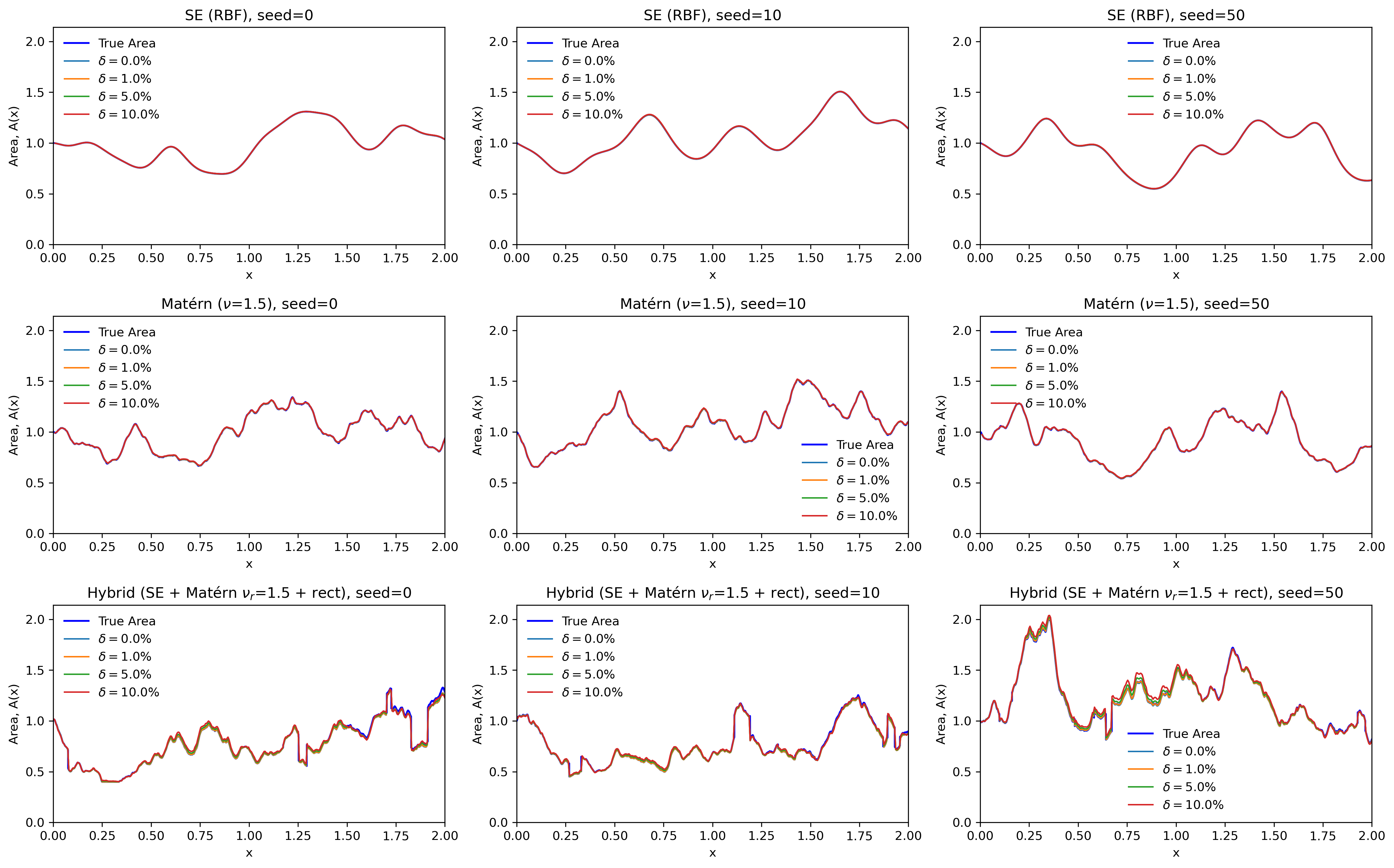}
\caption{\textbf{SG reconstructions under additive noise with non-reflecting boundary update.} Rows: smooth SE, rough Mat\'ern, hybrid. Columns: seeds $0,10,50$.}
\label{fig:sg_noise_examples_RB}
\end{figure}

Figure~\ref{fig:sg_noise_examples_RB} shows SG reconstructions at noise levels $\delta\in\{0,0.01,0.05,0.10\}$. Smooth and rough profiles remain closely aligned with the truth even at $\delta=10\%$. Hybrid profiles exhibit larger discrepancies as $\delta$ increases, reflecting the difficulty of resolving non-smooth features in noisy data. Table~\ref{tab:sg_noisy_rel_errors} shows that relative $L^2$-errors grow from $\sim10^{-4}$ (noiseless) to $\sim10^{-3}$ at $\delta=10\%$ for smooth profiles, while $H^1$-errors for hybrid profiles exceed $10^{-1}$.

\begin{figure}[htbp]
\centering
\includegraphics[width=\textwidth]{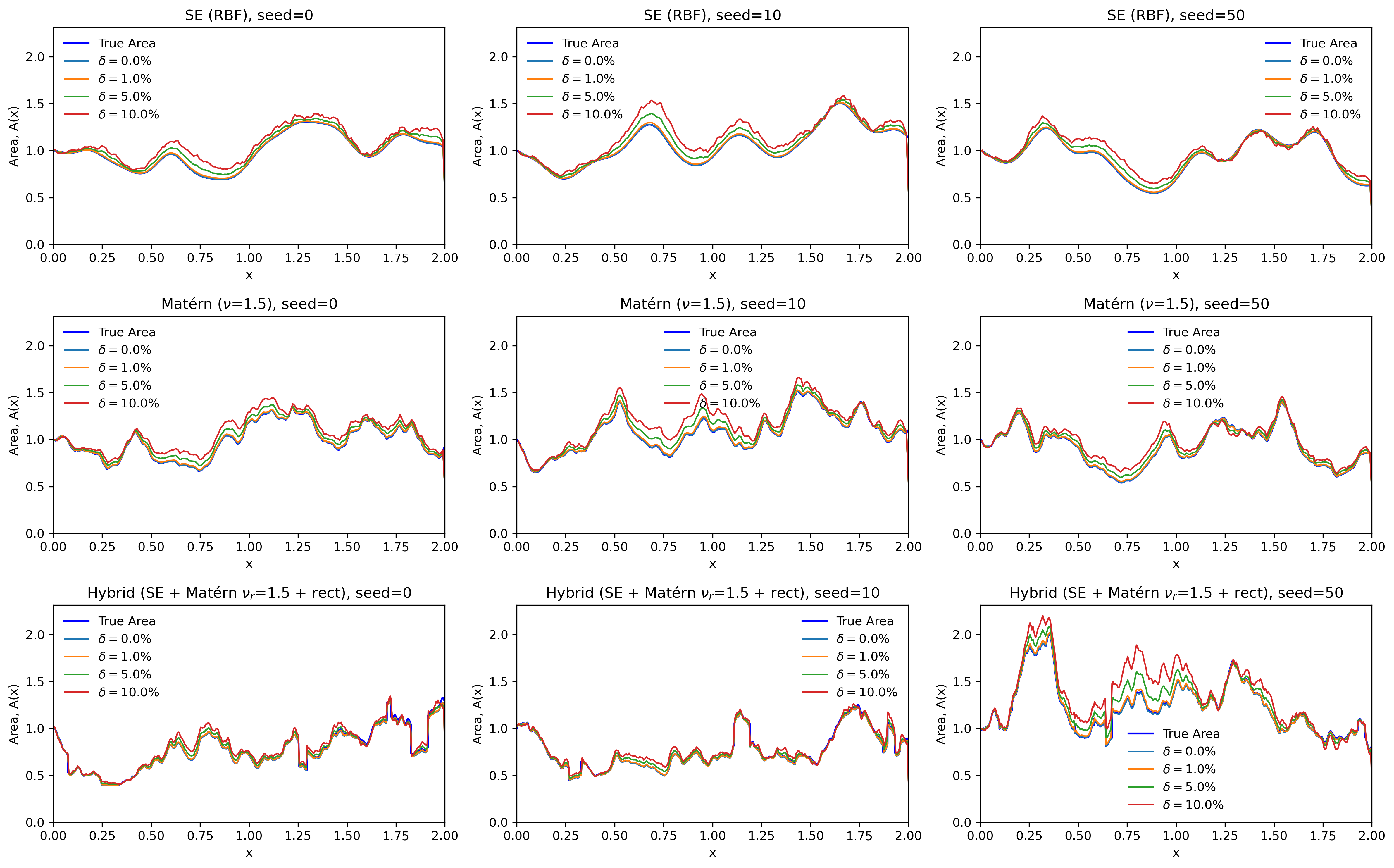}
\caption{\textbf{SG reconstructions without non-reflecting boundary update.} Boundary reflections inflate the noise norm, causing stronger oscillations at larger noise levels.}
\label{fig:sg_noise_examples_WRB}
\end{figure}

A key numerical issue is the coupling between right-boundary treatment and noise normalisation. Without the absorbing boundary update, spurious reflections inflate $\|\tilde{\mathbf{h}}\|_{2,\Delta t}$ and hence the absolute noise amplitude across the full data vector, leading to stronger oscillations in reconstructed profiles (Figure~\ref{fig:sg_noise_examples_WRB}).

\begin{table}[htbp]
\centering
\scriptsize
\setlength{\tabcolsep}{4pt}
\begin{tabular}{llcccccccc}
\toprule
\multirow{2}{*}{Seed} & \multirow{2}{*}{Case} & \multicolumn{2}{c}{$\delta=0\%$} & \multicolumn{2}{c}{$\delta=1\%$} & \multicolumn{2}{c}{$\delta=5\%$} & \multicolumn{2}{c}{$\delta=10\%$} \\
\cmidrule(lr){3-4} \cmidrule(lr){5-6} \cmidrule(lr){7-8} \cmidrule(lr){9-10}
& & $H^1_{\mathrm{rel}}$ & $L^2_{\mathrm{rel}}$ & $H^1_{\mathrm{rel}}$ & $L^2_{\mathrm{rel}}$ & $H^1_{\mathrm{rel}}$ & $L^2_{\mathrm{rel}}$ & $H^1_{\mathrm{rel}}$ & $L^2_{\mathrm{rel}}$ \\
\midrule
0 & SE (RBF) & $1.60\times 10^{-3}$ & $1.89\times 10^{-4}$ & $1.82\times 10^{-3}$ & $2.59\times 10^{-4}$ & $5.21\times 10^{-3}$ & $7.20\times 10^{-4}$ & $1.02\times 10^{-2}$ & $1.34\times 10^{-3}$ \\
  & Mat\'ern & $7.80\times 10^{-3}$ & $4.25\times 10^{-4}$ & $8.98\times 10^{-3}$ & $5.80\times 10^{-4}$ & $1.41\times 10^{-2}$ & $1.64\times 10^{-3}$ & $2.10\times 10^{-2}$ & $3.09\times 10^{-3}$ \\
  & Hybrid & $3.51\times 10^{-2}$ & $1.79\times 10^{-2}$ & $3.69\times 10^{-2}$ & $1.68\times 10^{-2}$ & $6.78\times 10^{-2}$ & $1.64\times 10^{-2}$ & $1.22\times 10^{-1}$ & $2.42\times 10^{-2}$ \\
\midrule
10 & SE (RBF) & $2.35\times 10^{-3}$ & $4.21\times 10^{-4}$ & $2.66\times 10^{-3}$ & $5.45\times 10^{-4}$ & $5.97\times 10^{-3}$ & $1.15\times 10^{-3}$ & $1.10\times 10^{-2}$ & $1.97\times 10^{-3}$ \\
  & Mat\'ern & $8.27\times 10^{-3}$ & $7.53\times 10^{-4}$ & $8.34\times 10^{-3}$ & $1.02\times 10^{-3}$ & $1.09\times 10^{-2}$ & $2.35\times 10^{-3}$ & $1.70\times 10^{-2}$ & $4.12\times 10^{-3}$ \\
  & Hybrid & $5.50\times 10^{-2}$ & $1.24\times 10^{-2}$ & $5.57\times 10^{-2}$ & $1.13\times 10^{-2}$ & $7.45\times 10^{-2}$ & $9.98\times 10^{-3}$ & $1.16\times 10^{-1}$ & $1.58\times 10^{-2}$ \\
\bottomrule
\end{tabular}
\caption{\textbf{Relative error analysis under additive noise for the SG algorithm.} Relative $H^1$- and $L^2$-errors for the smooth SE, rough Mat\'ern and hybrid profiles. The smooth and rough profiles remain stable under increasing noise, while the hybrid profiles exhibit larger errors, especially in the $H^1$-metric.}
\label{tab:sg_noisy_rel_errors}
\end{table}

\subsection{KLO algorithm}
\subsubsection{Noiseless case}
\begin{figure}[htbp]
\centering\includegraphics[width=\textwidth]{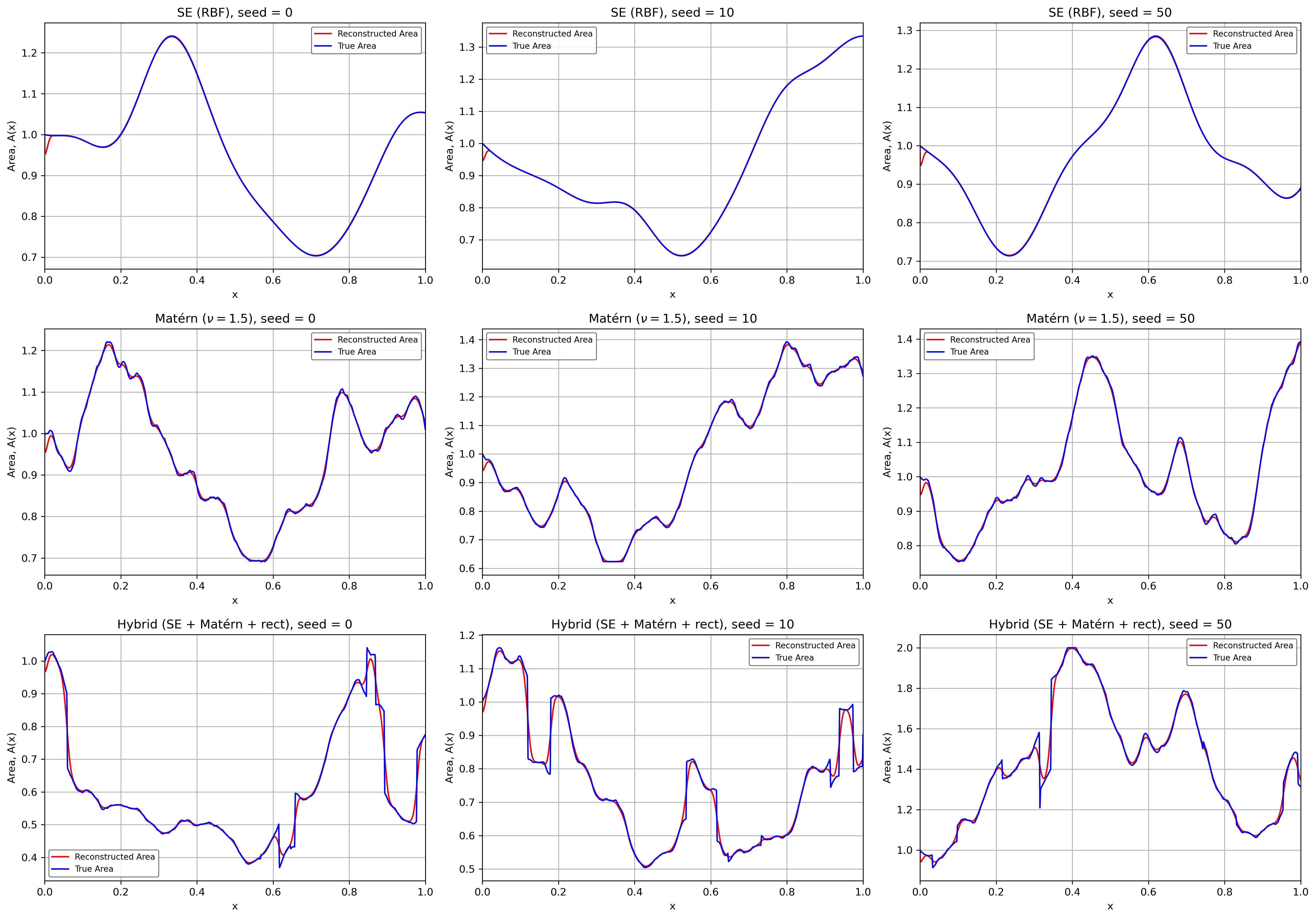}
\caption{\textbf{Uncorrected noise-free KLO reconstructions.}  The reconstructions recover the main interior features, but a  characteristic numerical artefact appears near the inlet $x=0$.}
\label{fig:bc_noiseless_examples}
\end{figure}

\begin{figure}[htbp]
    \centering
    \begin{subfigure}{\textwidth}
        \centering
        \includegraphics[width=\linewidth]{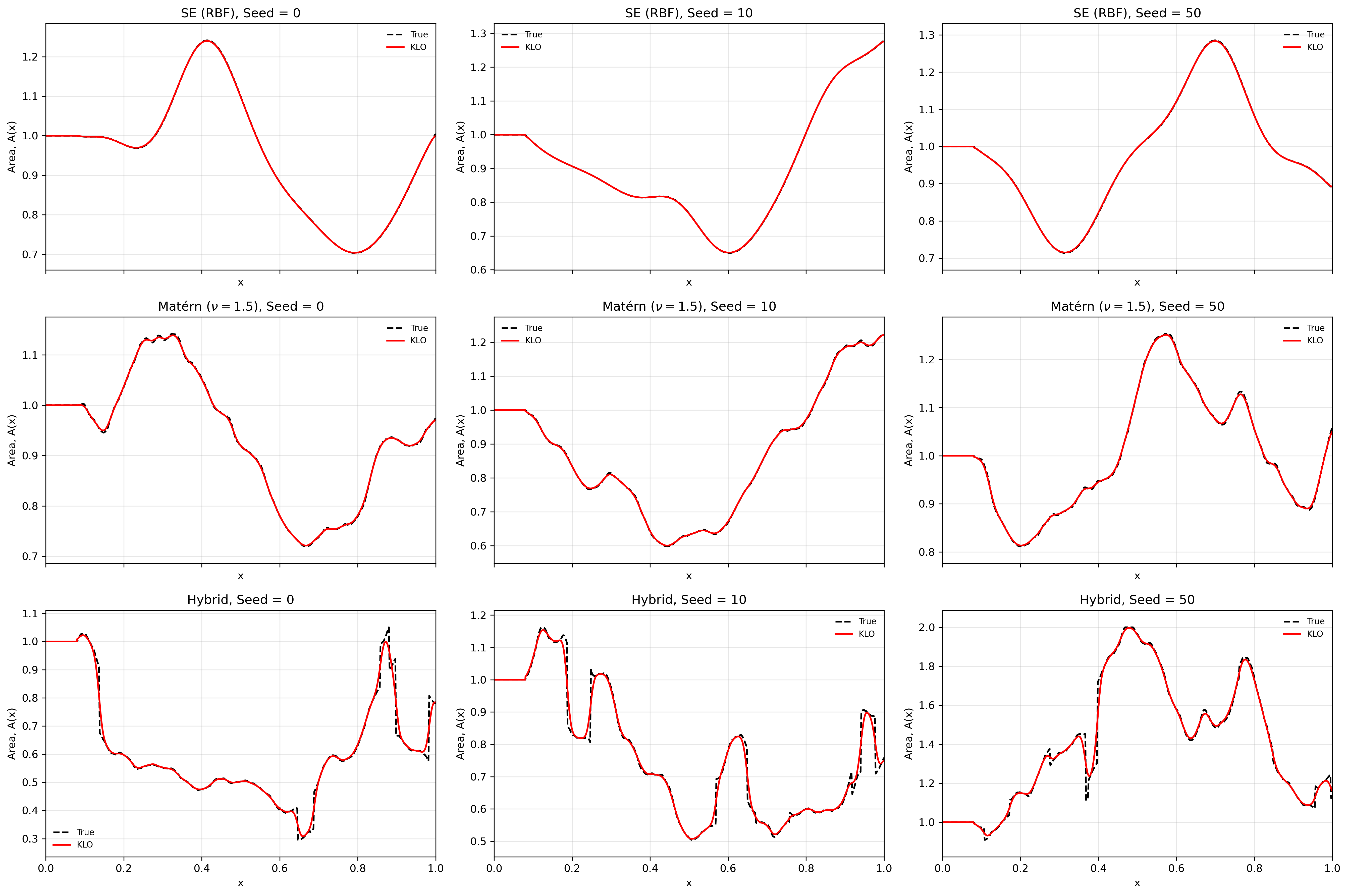}
        \caption{}
        \label{fig:klo_xin008}
    \end{subfigure}
    \begin{subfigure}{\textwidth}
        \centering
        \includegraphics[width=\linewidth]{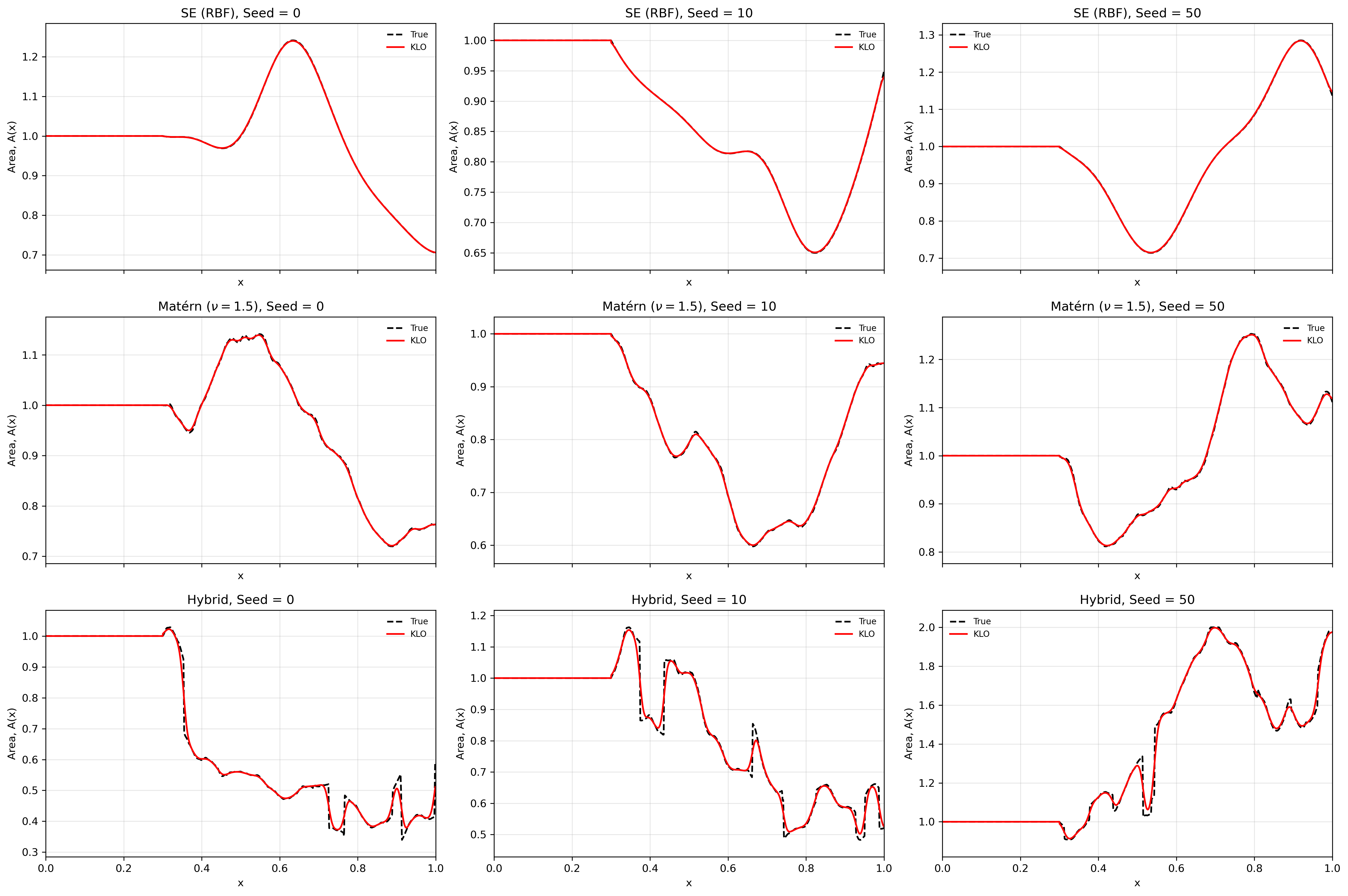}
        \caption{}
        \label{fig:klo_xin030}
    \end{subfigure}
    \caption{\textbf{Boundary corrected noise-free KLO reconstructions.} Shorter inlet correction removes the dominant artefact; longer correction imposes an artificial flat region.}
    \label{fig:bc_noiseless_examples_corr}
\end{figure}

Figure~\ref{fig:bc_noiseless_examples} shows noise-free KLO reconstructions. The method recovers interior geometry accurately for all three profile classes. However, a systematic numerical artefact appears near the inlet $x=0$ due to the small-$r$ projector being under-resolved; this is amplified in the $H^1$ norm (Table~\ref{tab:bc_noiseless_errors}). Imposing $A(x)=1$ for $x\leq x_{\min}$ (inlet regularisation with $x_{\min}=0.08$) removes the dominant artefact while preserving the main interior geometry; a larger value $x_{\min}=0.30$ over-regularises (Figure~\ref{fig:bc_noiseless_examples_corr}). Dense-grid evaluations reveal this boundary-layer structure more clearly than coarse grids, which can mask the defect and produce artificially optimistic $H^1$-error estimates.

\begin{table}[htbp]
\centering
\scriptsize
\setlength{\tabcolsep}{6pt}
\begin{tabular}{llcccc}
\toprule
Seed & Case & $H^1_{\mathrm{abs}}$ & $H^1_{\mathrm{rel}}$ & $L^2_{\mathrm{abs}}$ & $L^2_{\mathrm{rel}}$ \\
\midrule
0 & SE (RBF) & $3.7793\times 10^{-1}$ & $2.1839\times 10^{-1}$ & $3.9549\times 10^{-3}$ & $4.0969\times 10^{-3}$ \\
  & Matérn & $4.9192\times 10^{-1}$ & $2.9871\times 10^{-1}$ & $4.2821\times 10^{-3}$ & $4.5336\times 10^{-3}$ \\
  & Hybrid & $1.0177\times 10^{1}$ & $9.0846\times 10^{-1}$ & $1.9888\times 10^{-2}$ & $3.1482\times 10^{-2}$ \\
\midrule
10 & SE (RBF) & $4.1990\times 10^{-1}$ & $2.6333\times 10^{-1}$ & $4.2109\times 10^{-3}$ & $4.4202\times 10^{-3}$ \\
   & Matérn & $5.3327\times 10^{-1}$ & $3.0032\times 10^{-1}$ & $4.4756\times 10^{-3}$ & $4.8990\times 10^{-3}$ \\
   & Hybrid & $1.0548\times 10^{1}$ & $9.1042\times 10^{-1}$ & $2.0217\times 10^{-2}$ & $2.5939\times 10^{-2}$ \\
\bottomrule
\end{tabular}
\caption{\textbf{Baseline errors for uncorrected noise-free KLO reconstructions.} The elevated \(H^1\)-errors are caused by the inlet artefact.}
\label{tab:bc_noiseless_errors}
\end{table}

\subsubsection{Performance in the noisy case} \label{Noise_4_KLO_results}

\begin{figure}[htbp]
    \centering
    \includegraphics[width=\textwidth]{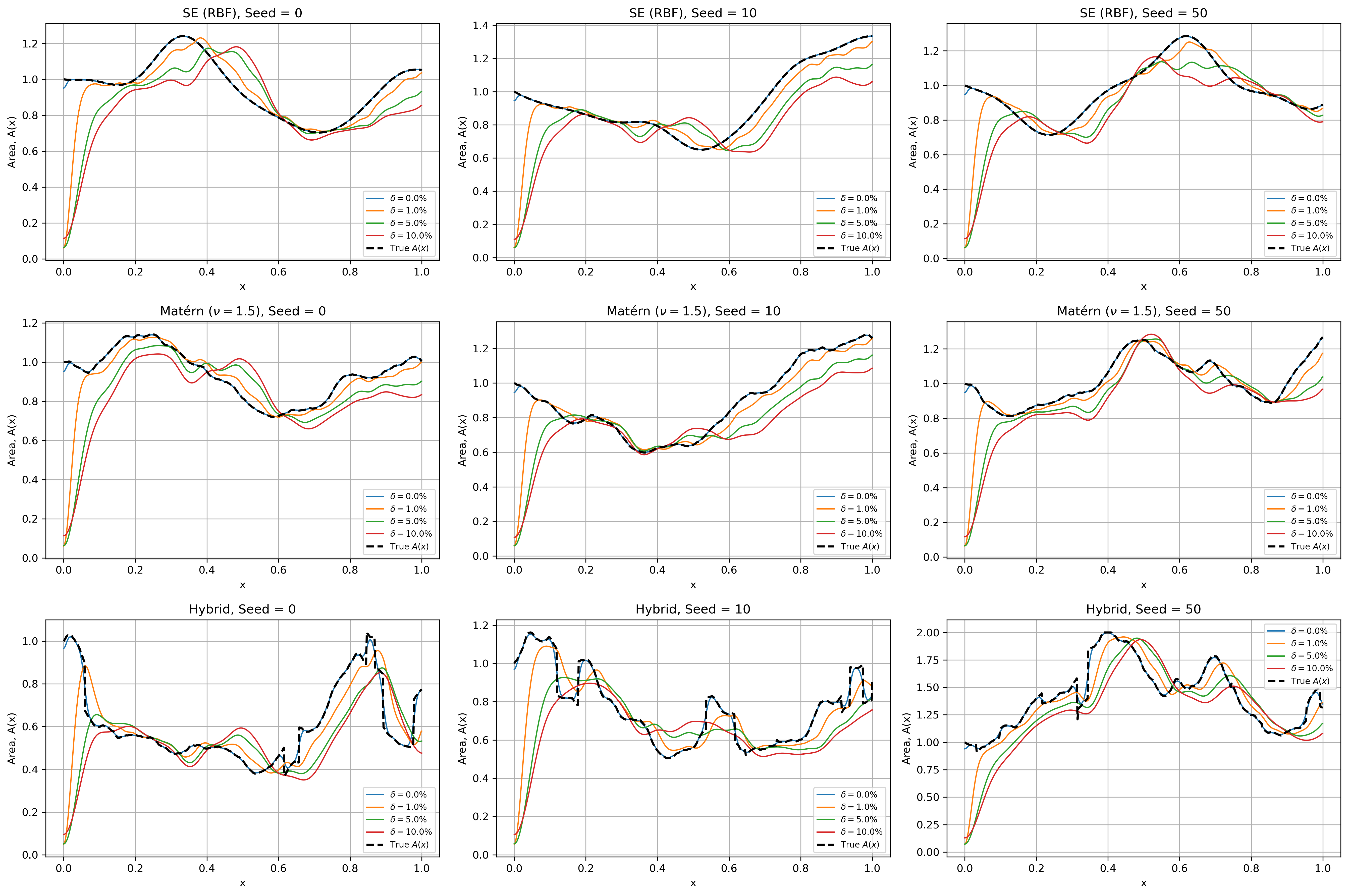}
    \caption{KLO reconstructions under additive noise ($\delta\in\{0\%,1\%,5\%,10\%\}$) for smooth SE, rough Mat\'ern, and hybrid profiles.}
    \label{fig:bc_noisy_examples}
\end{figure}

Figure~\ref{fig:bc_noisy_examples} shows KLO reconstructions under additive noise. The inlet boundary layer becomes more pronounced as $\delta$ increases. Away from the inlet, large-scale geometric structure is retained for SE and Mat\'ern profiles at low noise, but degrades rapidly for hybrid profiles. Table~\ref{tab:bc_noisy_rel_errors} shows $L^2$-errors growing monotonically with $\delta$, while $H^1$-errors are non-monotone: after a sharp increase at $\delta=1\%$ the stronger regularisation and smoothing at larger $\delta$ suppress high-frequency oscillations, reducing the $H^1$-error even as the $L^2$-error continues to grow.

\begin{table}[htbp]
\centering
\scriptsize
\setlength{\tabcolsep}{4pt}
\begin{tabular}{llcccccccc}
\toprule
\multirow{2}{*}{\textbf{Seed}} & \multirow{2}{*}{\textbf{Case}} & \multicolumn{2}{c}{$\delta=0\%$} & \multicolumn{2}{c}{$\delta=1\%$} & \multicolumn{2}{c}{$\delta=5\%$} & \multicolumn{2}{c}{$\delta=10\%$} \\
\cmidrule(lr){3-4} \cmidrule(lr){5-6} \cmidrule(lr){7-8} \cmidrule(lr){9-10}
& & $H^1_{\mathrm{rel}}$ & $L^2_{\mathrm{rel}}$ & $H^1_{\mathrm{rel}}$ & $L^2_{\mathrm{rel}}$ & $H^1_{\mathrm{rel}}$ & $L^2_{\mathrm{rel}}$ & $H^1_{\mathrm{rel}}$ & $L^2_{\mathrm{rel}}$ \\
\midrule
0 & SE (RBF) & $2.18\times 10^{-1}$ & $4.10\times 10^{-3}$ & $2.39\times 10^{0}$ & $1.43\times 10^{-1}$ & $1.81\times 10^{0}$ & $2.22\times 10^{-1}$ & $1.68\times 10^{0}$ & $2.52\times 10^{-1}$ \\
  & Mat\'ern & $2.99\times 10^{-1}$ & $4.53\times 10^{-3}$ & $2.60\times 10^{0}$ & $1.43\times 10^{-1}$ & $1.87\times 10^{0}$ & $2.15\times 10^{-1}$ & $1.63\times 10^{0}$ & $2.43\times 10^{-1}$ \\
  & Hybrid   & $9.08\times 10^{-1}$ & $3.15\times 10^{-2}$ & $1.07\times 10^{0}$ & $2.60\times 10^{-1}$ & $1.06\times 10^{0}$ & $3.49\times 10^{-1}$ & $1.04\times 10^{0}$ & $3.66\times 10^{-1}$ \\
\midrule
10 & SE (RBF) & $2.63\times 10^{-1}$ & $4.42\times 10^{-3}$ & $2.63\times 10^{0}$ & $1.43\times 10^{-1}$ & $1.95\times 10^{0}$ & $2.22\times 10^{-1}$ & $1.77\times 10^{0}$ & $2.56\times 10^{-1}$ \\
   & Mat\'ern & $3.00\times 10^{-1}$ & $4.90\times 10^{-3}$ & $2.44\times 10^{0}$ & $1.49\times 10^{-1}$ & $1.82\times 10^{0}$ & $2.23\times 10^{-1}$ & $1.59\times 10^{0}$ & $2.53\times 10^{-1}$ \\
   & Hybrid   & $9.10\times 10^{-1}$ & $2.59\times 10^{-2}$ & $1.04\times 10^{0}$ & $2.12\times 10^{-1}$ & $1.02\times 10^{0}$ & $3.07\times 10^{-1}$ & $1.02\times 10^{0}$ & $3.31\times 10^{-1}$ \\
\bottomrule
\end{tabular}
\caption{\textbf{Relative $H^1$- and $L^2$-errors for noisy KLO reconstructions (Seeds $0$ and $10$).}} 
\label{tab:bc_noisy_rel_errors}
\end{table}

\subsection{Cross algorithm validation}

\subsubsection{Reconstruction from SG-generated data}
\begin{figure}[htpb]
\centering
\includegraphics[width=\linewidth]{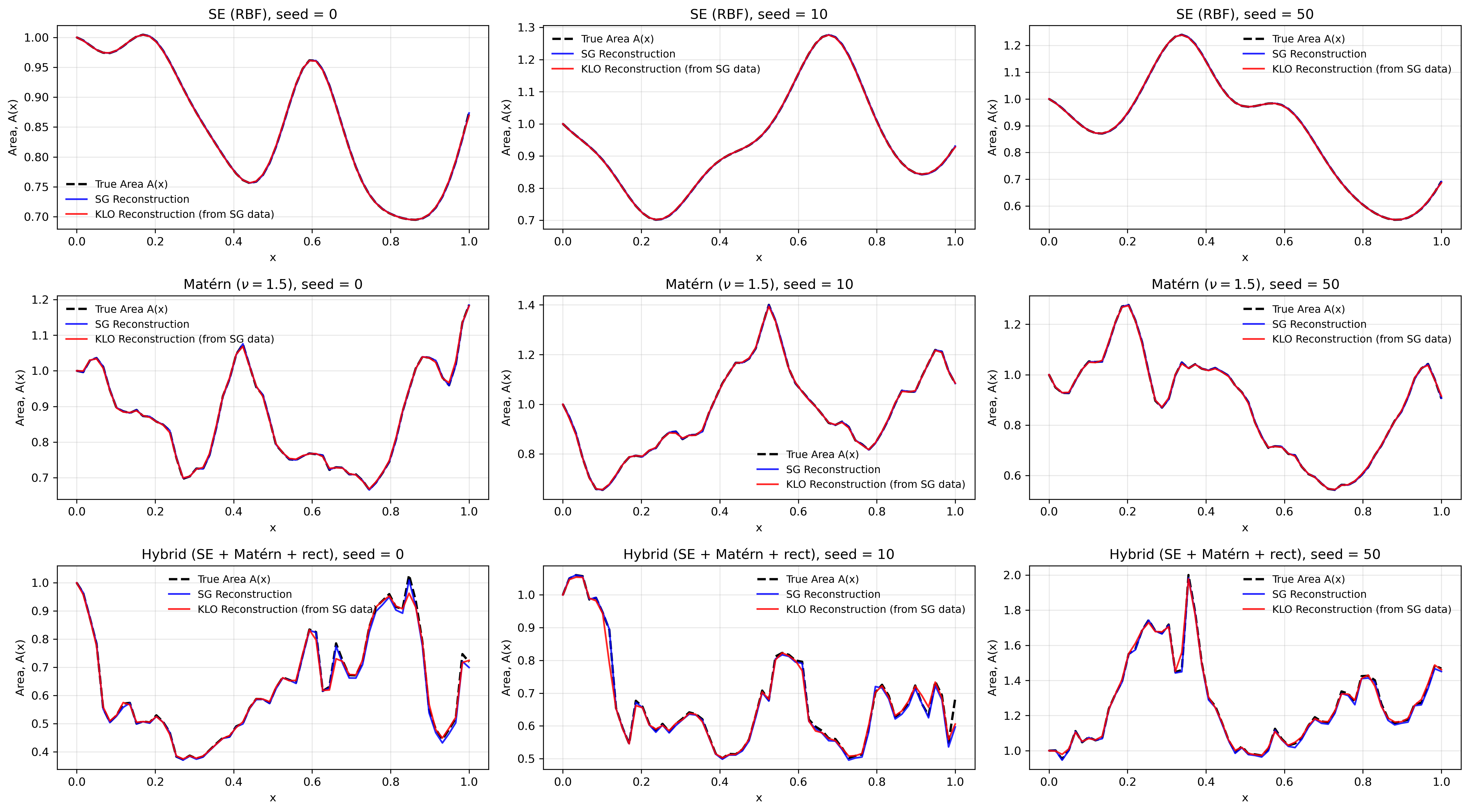}
\caption{\textbf{Cross-validation: KLO reconstruction (red) from SG-generated data, compared with true area profile (dashed black) and the native SG reconstruction (blue).}}
\label{SG_data_in_KLO}
\end{figure}

Figure~\ref{SG_data_in_KLO} shows KLO reconstructions from SG-generated data. For SE and Mat\'ern profiles, KLO remains in close agreement with both the true profile and the native SG reconstruction. Hybrid profiles show small local discrepancies near sharp transitions but preserve the main geometric structure.

\subsubsection{Reconstruction from KLO-generated data}
\begin{figure}[htpb]
\centering
\includegraphics[width=\linewidth]{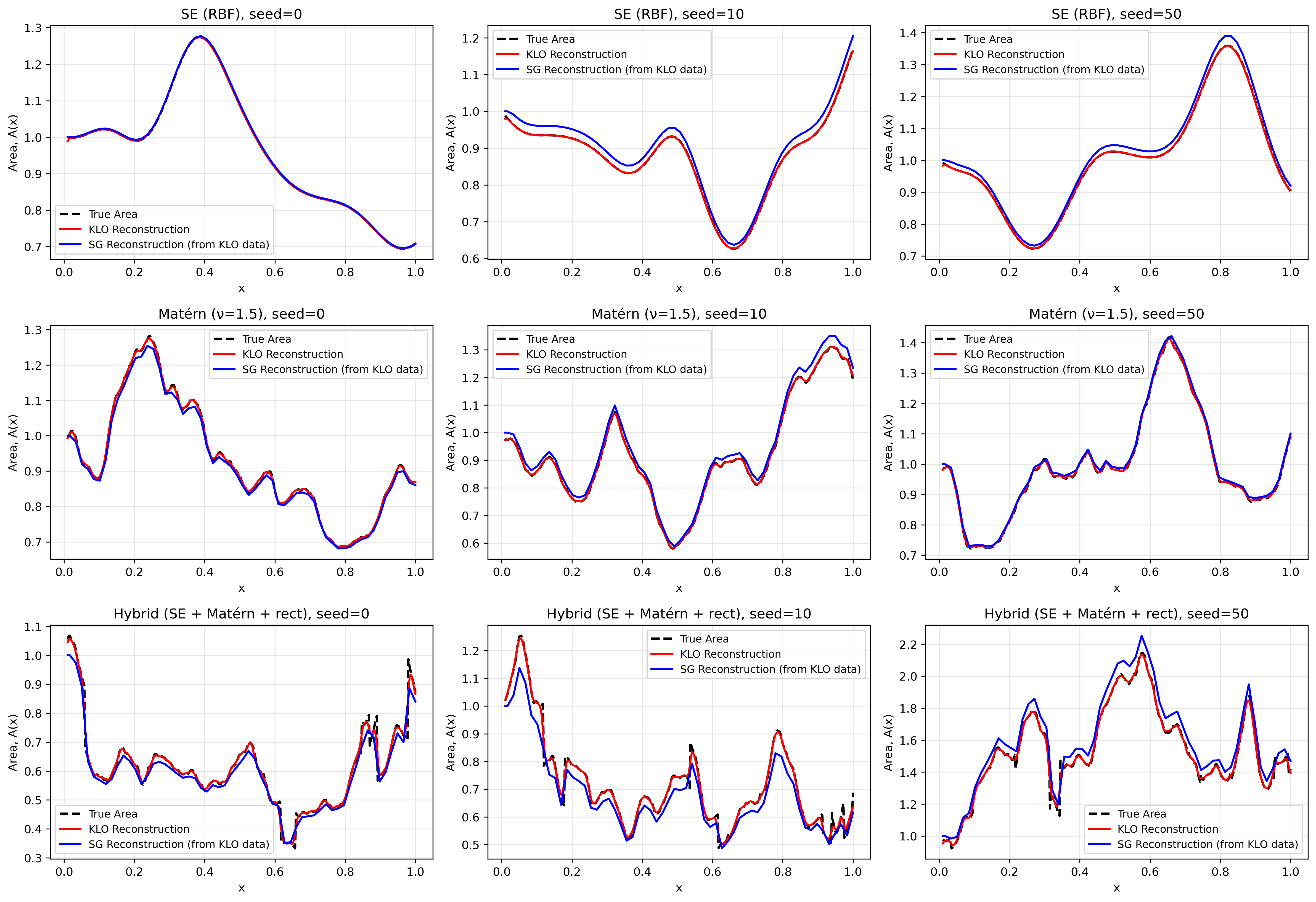}
\caption{\textbf{Cross-validation: SG reconstruction (blue) from KLO-generated data, compared with true profile (dashed black)} and native KLO reconstruction (red).}
\label{KLO_data_in_SG}
\end{figure}

Figure~\ref{KLO_data_in_SG} shows SG reconstructions from KLO-generated data. SE and Mat\'ern profiles are recovered consistently. Hybrid profiles exhibit larger local deviations, attributable to the required temporal differentiation step \eqref{hSG_from_KLO}, which amplifies numerical irregularities. Overall, both methods remain broadly consistent when driven by the other's data.

\subsubsection{Statistical validation in the SG-to-KLO cross-transfer setting}\label{sec:statistical-validation}
To complement the representative examples and assess whether the observed trends persist across a larger ensemble, we perform a statistical comparison of the SG and KLO reconstructions in the SG-to-KLO transfer setting. This transfer direction was selected because the conversion it requires is numerically stable, whereas the reverse direction requires differentiation of the boundary trace and is correspondingly more sensitive to perturbation. See the discussion at the end of Subsection~\ref{sec:klo-from-sg}. We consider an ensemble of \(n=1000\) randomly generated area profiles and compare the resulting SG and KLO reconstructions at prescribed noise levels. 

\begin{figure}[htbp]
\centering
\includegraphics[width=\textwidth]{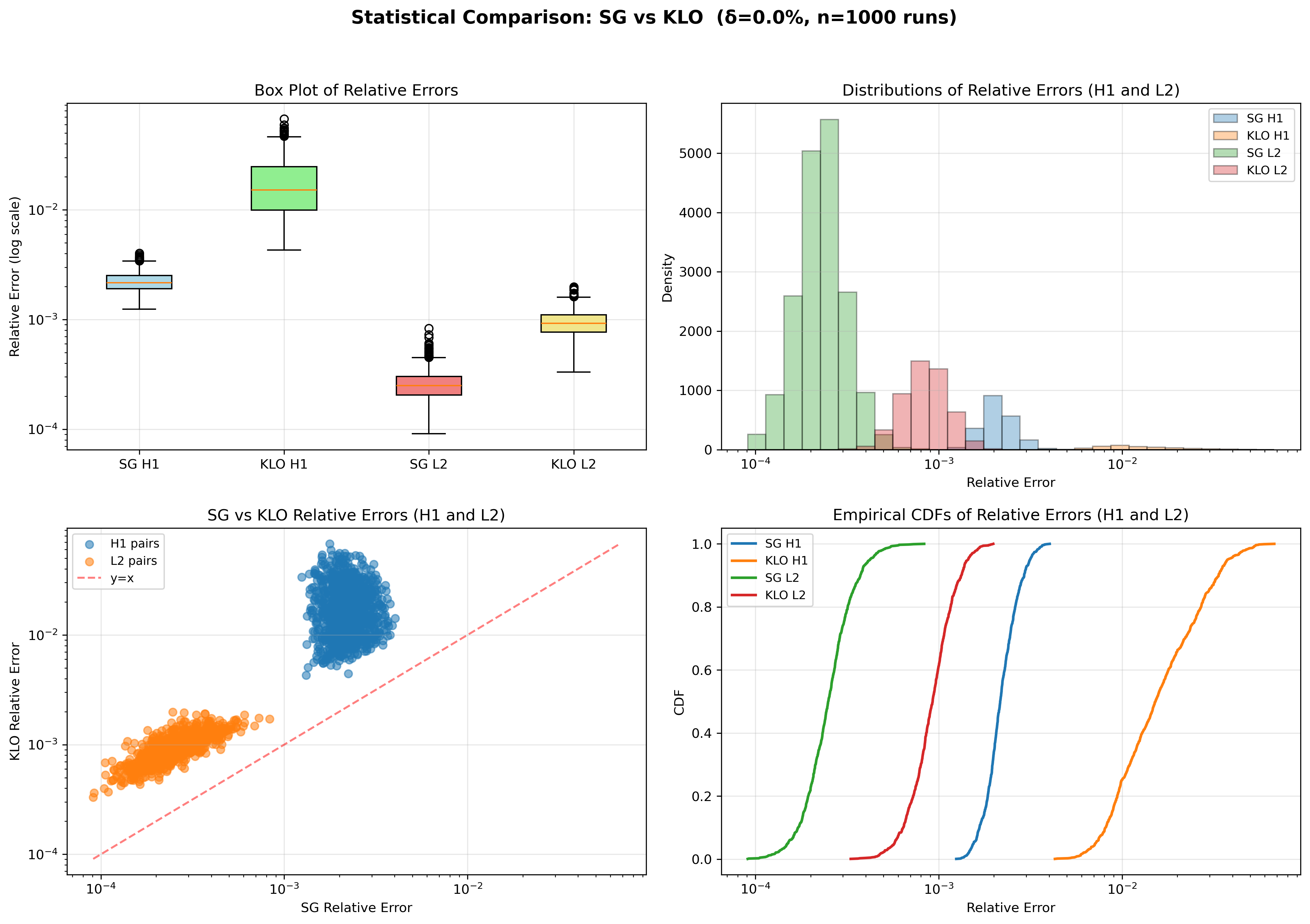}
\caption{\textbf{Statistical comparison of SG and KLO in the noiseless case (\(\delta=0\%, n=1000\)). Distributional analysis.} The panels show relative error box plots, histogram distributions, empirical cumulative distribution functions, and paired scatter plots for the \(H^1\)- and \(L^2\)-metrics. The separation of the distributions and the location of the paired samples above the reference line \(y=x\) indicate consistently smaller SG errors.}
\label{fig:comparison_noiseless1}
\end{figure}

\paragraph{Noiseless case ($\delta=0\%$).}
As shown in Figure~\ref{fig:comparison_noiseless1}, SG consistently yields smaller errors than KLO across all $1000$ realisations. The error distributions are shifted toward lower values for SG, and the paired scatter plots show all observations above the reference line, indicating that every realisation favours SG. The mean paired differences are $\bar{d}_{H^1}=-1.61\times10^{-2}$ (95\% CI: $[-1.68,-1.55]\times10^{-2}$) and $\bar{d}_{L^2}=-6.89\times10^{-4}$ (95\% CI: $[-7.03,-6.76]\times10^{-4}$). Both paired $t$-tests and Wilcoxon tests give $p<0.001$, Cohen's $d=-1.54$ for $H^1$ and $-3.26$ for $L^2$, and the rank-biserial correlation is $-1.000$.

\begin{figure}[htbp]
\centering
\includegraphics[width=\textwidth]{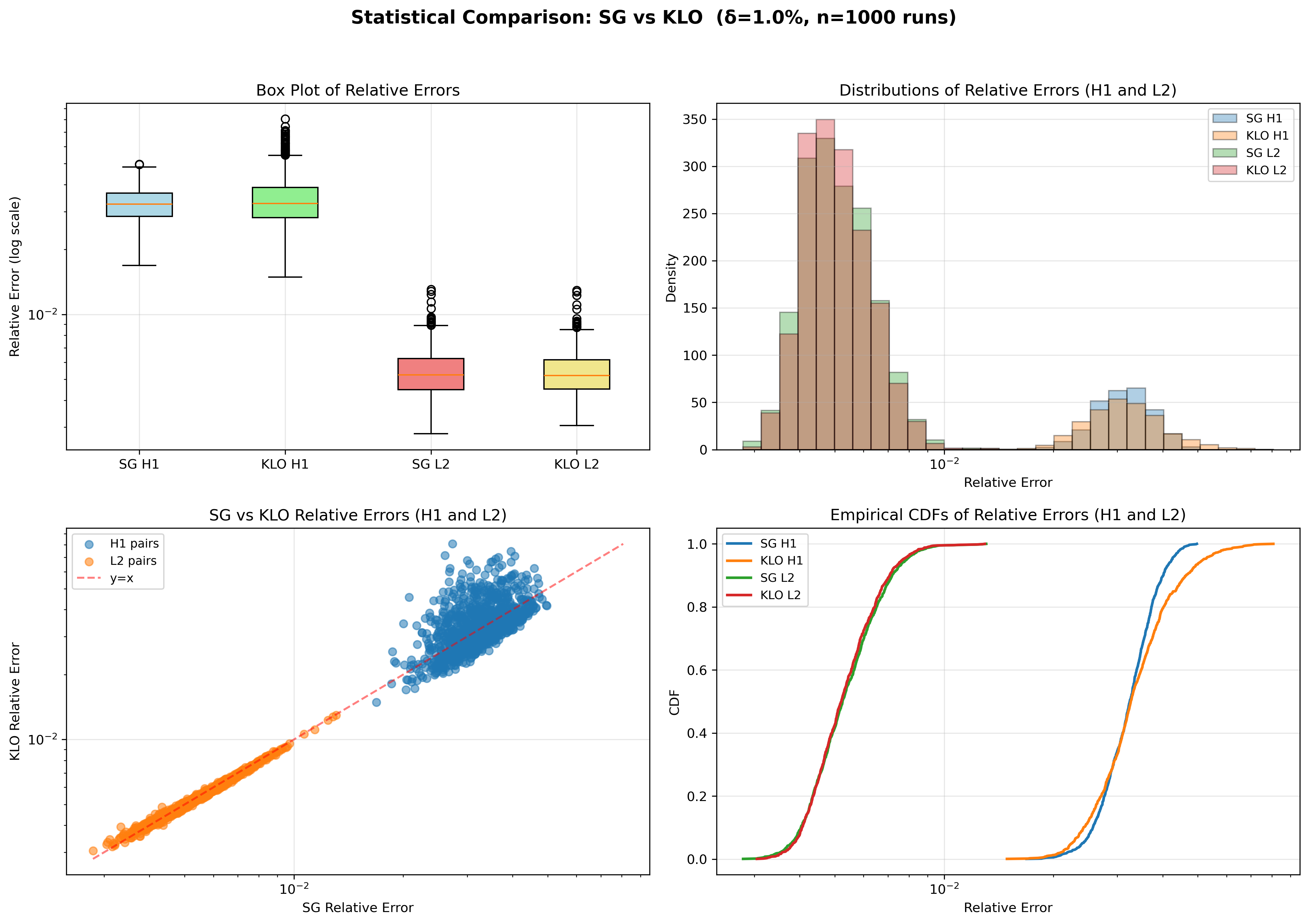}

\caption{\textbf{Statistical comparison under \(1\%\) noise (\(n=1000\)). Distributional analysis.} Box plots, histograms, empirical cumulative distribution functions, and paired scatter plots show substantial overlap, with a heavier upper tail for KLO, between SG and KLO in the \(H^1\)-metric, while the \(L^2\)-metric begins to favour KLO.}
\label{fig:comparison_1pct1}
\end{figure}

\paragraph{Low noise case \(\boldsymbol{(\delta=1\%)}\).} Figure~\ref{fig:comparison_1pct1} shows that, at $1\%$ noise, the clear noiseless ordering disappears for the $H^1$-metric: the error distributions overlap and the paired scatter points fall on both sides of the reference line. The two summaries of the paired difference disagree in sign: the mean is $\bar{d}_{H^1}=-1.76\times10^{-3}$ (95\% CI: $[-2.26,-1.26]\times10^{-3}$), favouring SG, while the median is $8.78\times10^{-4}$ (bootstrap 95\% CI: $[3.59,13.11]\times10^{-4}$), favouring KLO, and KLO attains the smaller $H^1$-error in $55.8\%$ of realisations. The Wilcoxon signed-rank test does not reject the null hypothesis ($p=0.119$), the only such case across the noise levels examined, and the rank-biserial correlation is $-0.057$. The paired $t$-test is significant ($p<0.001$) but the effect size is small (Cohen's $d=-0.219$), reflecting the large sample rather than a practically meaningful difference. This pattern indicates that KLO is typically marginally more accurate but occasionally incurs large errors that dominate the mean. For $L^2$, the figures show a clearer advantage for KLO: $\bar{d}_{L^2}=5.60\times10^{-5}$, $p<0.001$, KLO win rate $67.6\%$, and Cohen's $d=0.445$.

\begin{figure}[htbp]
\centering
\includegraphics[width=\textwidth]{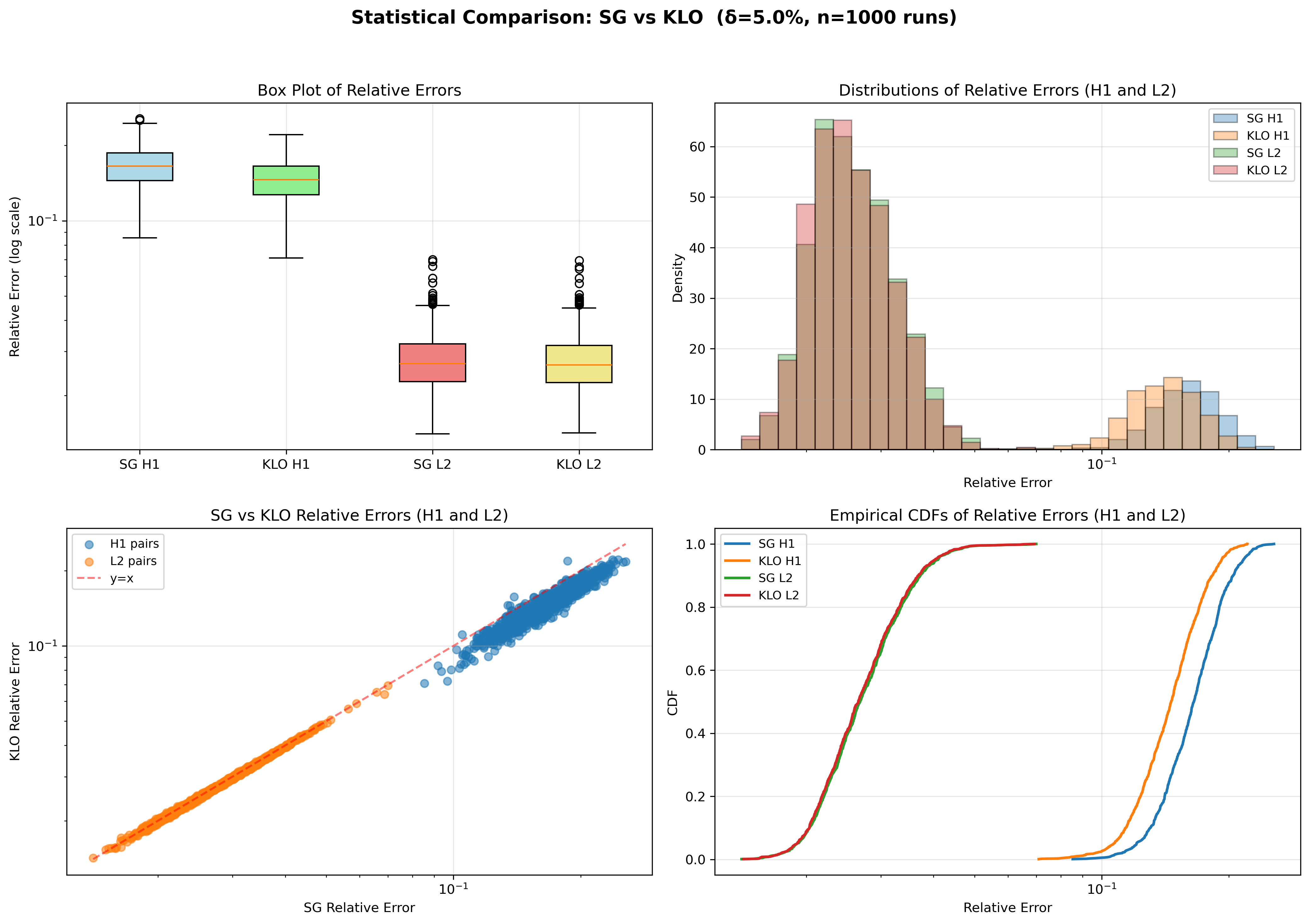}
\caption{\textbf{Statistical comparison at \(5\%\) noise (\(n=1000\)). Distributional analysis.} Box plots, histograms, empirical cumulative distribution functions, and paired scatter plots compare the relative \(H^1\)- and \(L^2\)-errors of SG and KLO. The \(H^1\)-errors show a clear shift in favour of KLO, while the \(L^2\)-errors display a smaller but consistent improvement.}
\label{fig:comparison_5pct1}
\end{figure}

\paragraph{Moderate-noise case \(\boldsymbol{(\delta=5\%)}\).}Figure~\ref{fig:comparison_5pct1} shows a clear shift in favour of KLO at $5\%$ noise, especially in the $H^1$-metric: most paired scatter points lie below the reference line, indicating lower KLO errors. The mean paired difference is $\bar{d}_{H^1}=1.93\times10^{-2}$, with $p<0.001$, Cohen's $d=2.23$, and KLO win rate $99.0\%$. The $L^2$ results show a smaller but still robust advantage for KLO: $\bar{d}_{L^2}=3.11\times10^{-4}$, $p<0.001$, Cohen's $d=0.959$, and KLO win rate $84.9\%$.

\begin{figure}[htbp]
\centering
\includegraphics[width=\textwidth]{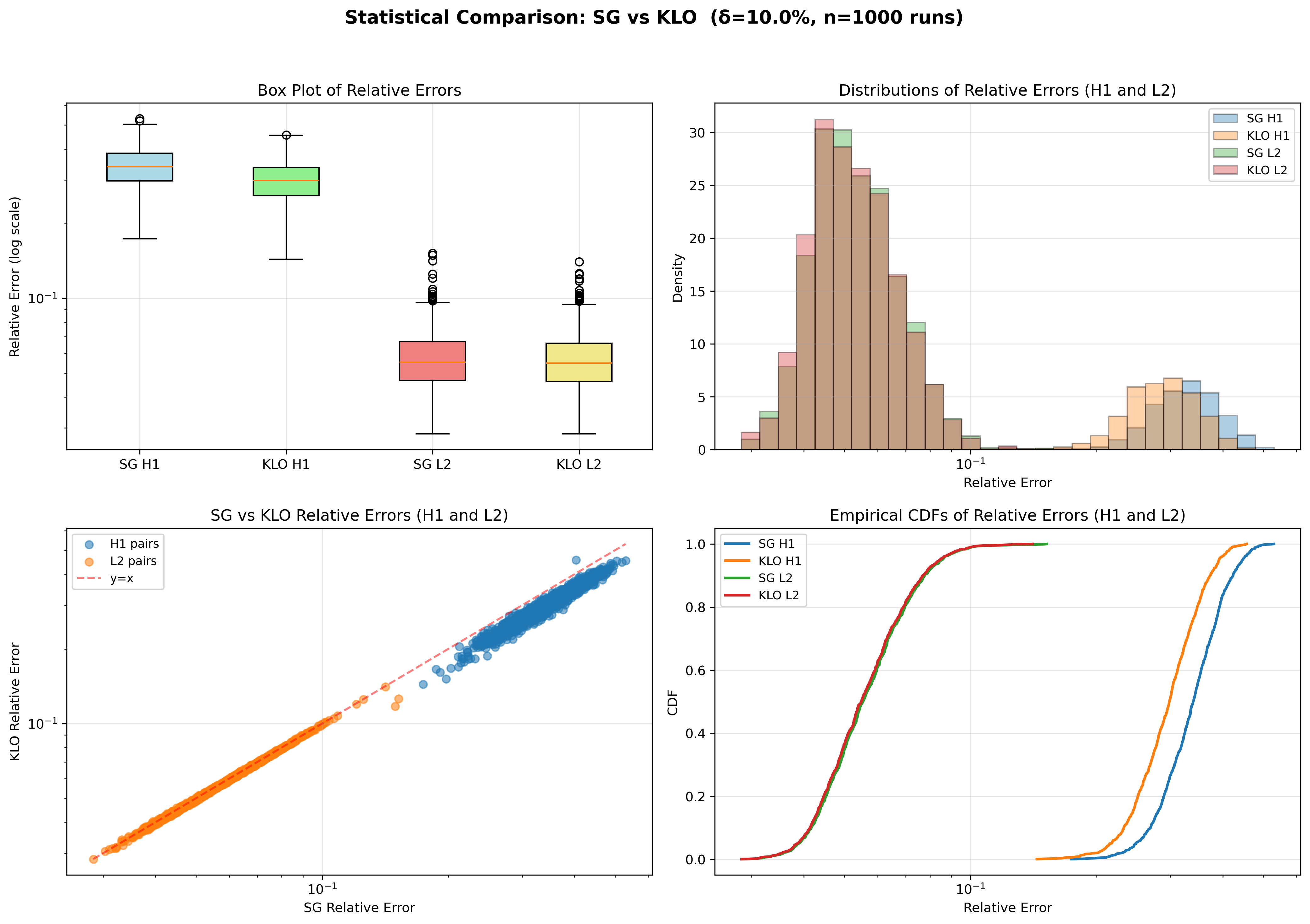}
\caption{\textbf{Statistical comparison at \(10\%\) noise (\(n=1000\)). Distributional analysis.} Box plots, histograms, empirical cumulative distribution functions, and paired scatter plots compare SG and KLO relative errors. KLO remains favoured in the \(H^1\)-metric, while the \(L^2\)-metric shows a smaller but consistent shift toward lower KLO errors.}
\label{fig:comparison_10pct1}
\end{figure}

\paragraph{High noise case \(\boldsymbol{(\delta=10\%)}\).} Figure~\ref{fig:comparison_10pct1} shows that the advantage of KLO observed at $5\%$ noise persists and becomes stronger at $10\%$ noise. For $H^1$, the KLO error distribution is shifted toward lower values and nearly all paired points favour KLO: $\bar{d}_{H^1}=4.23\times10^{-2}$, $p<0.001$, Cohen's $d=3.05$, and KLO win rate $99.9\%$. For $L^2$, KLO also remains favoured, with $\bar{d}_{L^2}=6.64\times10^{-4}$, $p<0.001$, Cohen's $d=0.490$, and KLO win rate $90.1\%$.

The systematic evolution of these performance metrics across the full statistical ensemble is summarised in Table~\ref{tab:stat_summary}, which highlights the noise-dependent crossover in comparative performance.

\begin{table}[htbp]
\centering
\begin{tabular}{llcccc}
\toprule
$\delta$ & Metric & Mean (SG) & Mean (KLO) & Ratio (KLO/SG) & KLO Win Rate (\%) \\
\midrule
0\% & $H^1_{rel}$ & $2.24 \times 10^{-3}$ & $1.83 \times 10^{-2}$ & 8.19 & 0.0\% \\
 & $L^2_{rel}$ & $2.63 \times 10^{-4}$ & $9.52 \times 10^{-4}$ & 3.62 & 0.0\% \\
\midrule
1\% & $H^1_{rel}$ & $3.26 \times 10^{-2}$ & $3.44 \times 10^{-2}$ & 1.05 & 55.8\% \\
 & $L^2_{rel}$ & $5.48 \times 10^{-3}$ & $5.43 \times 10^{-3}$ & 0.99 & 67.6\% \\
\midrule
5\% & $H^1_{rel}$ & $1.66 \times 10^{-1}$ & $1.47 \times 10^{-1}$ & 0.88 & 99.0\% \\
 & $L^2_{rel}$ & $2.81 \times 10^{-2}$ & $2.77 \times 10^{-2}$ & 0.99 & 84.9\% \\
\midrule
10\% & $H^1_{rel}$ & $3.42 \times 10^{-1}$ & $3.00 \times 10^{-1}$ & 0.88 & 99.9\% \\
 & $L^2_{rel}$ & $5.81 \times 10^{-2}$ & $5.75 \times 10^{-2}$ & 0.99 & 90.1\% \\
\bottomrule
\end{tabular}
\caption{\textbf{Summary of SG and KLO statistical performance across noise levels (\(n=1000\), SG-to-KLO transfer setting).}  Ratio $<1$ indicates lower mean KLO error.}
\label{tab:stat_summary}
\end{table}

\subsection{Influence of the forward discretisation}
\label{sec:forward-convergence}

The reconstructions compared in Subsection 4.3.3 are driven by a single simulated boundary trace, from which the KLO Neumann-to-Dirichlet map is assembled as described in Subsection 3.3.2. The discretisation error of that forward solve therefore enters both reconstructions. In the noise-free case, it is the only error source apart from regularisation, so the errors reported at $\delta = 0$ are limited by the resolution rather than by the inversion methods themselves.

To quantify this, we refine the spatial and temporal grids
together and repeat the noise-free reconstruction. The measurement step
$\Delta t$ is halved twice, the forward grid is refined with it at a fixed
ratio of four, and the $r$-grid on which $s(\alpha, r)$ is sampled is refined so that $\Delta r / \Delta t$ is preserved. The smoothing length applied to $k(r)$ is held fixed in physical units. The finest level coincides with the configuration used in Subsection 4.3.3. Table~\ref{tab:convergence} reports the mean relative errors over one hundred randomly generated squared-exponential profiles.
\begin{table}[htbp]
\centering

\begin{tabular}{lccccc}
\toprule
$\Delta t$ & $N_t$ & SG $H^1_{\mathrm{rel}}$ & KLO $H^1_{\mathrm{rel}}$
& SG $L^2_{\mathrm{rel}}$ & KLO $L^2_{\mathrm{rel}}$ \\
\midrule
$1.0\times10^{-3}$ & $2050$ & $8.77\times10^{-3}$ & $5.04\times10^{-2}$
& $1.07\times10^{-3}$ & $4.40\times10^{-3}$ \\
$5.0\times10^{-4}$ & $4100$ & $4.37\times10^{-3}$ & $2.81\times10^{-2}$
& $5.33\times10^{-4}$ & $2.08\times10^{-3}$ \\
$2.5\times10^{-4}$ & $8200$ & $2.20\times10^{-3}$ & $1.80\times10^{-2}$
& $2.67\times10^{-4}$ & $9.70\times10^{-4}$ \\
\bottomrule
\end{tabular}
\caption{\textbf{Refinement study of the noise-free
reconstructions in the SG-to-KLO transfer setting.} Mean relative errors over one hundred profiles. The finest resolution is the one used in Subsection 4.3.3.}
\label{tab:convergence}
\end{table}

Both methods converge under refinement. For SG the error halves at each step, with observed ratios between $1.92$ and $2.08$ and median $2.00$ in both metrics, corresponding to first-order convergence. For KLO the relative $L^2$ error behaves similarly, with median ratio $2.14$, while the relative $H^1$ error converges more slowly, with ratios between $1.33$ and $2.51$ and median $1.76$. In no case does the error stop decreasing. The noise-free errors are thus a resolvable numerical floor, and since the scheme is consistent and the error behaves accordingly under refinement, a discretisation of comparable order would give the same limit.

Two further points explain why this floor does not affect the comparison between the methods. First, the two algorithms are supplied with identical simulated data, so the forward error is common to them and cannot favour either. The SG reconstruction is the more accurate of the two in all three hundred runs of Table~\ref{tab:convergence}, and refining the grid by a factor of four does not alter this. Second, the comparative conclusions of Subsection 4.3.3 concern behaviour under perturbation, where the reconstruction errors are far larger than this floor. Taking the SG mean relative $H^1$ error from Table 6, the noise-free value is $2.24\times10^{-3}$, against $3.26\times10^{-2}$, $1.66\times10^{-1}$ and $3.42\times10^{-1}$ at $\delta = 1\%$, $5\%$ and $10\%$, that is, larger by factors of $14.6$, $74.1$ and $152.7$. The noise-dependent crossover therefore occurs well above the level at which the forward discretisation is significant.

\section{Conclusion}\label{section5:conclusion}

The SG algorithm is robust and accurate for geometric reconstruction. Under noise-free conditions it achieves high precision for smooth and rough profiles and establishes a strong benchmark. It remains stable under additive noise up to $10\%$, recovering the dominant structure of all three profile classes. The main stability concern is the coupling between the right-boundary treatment and the relative-noise normalisation: reflections without the absorbing update inflate the effective noise amplitude, so the non-reflecting boundary condition is essential for reliable noisy performance. Hybrid profiles are consistently the most demanding due to their sharp local transitions.

The KLO method recovers interior geometry accurately but is affected by a persistent inlet boundary effect associated with the small-$r$ projector. Inlet regularisation (prescribing $A(x)=1$ near $x=0$) substantially reduces this artefact. Dense-grid evaluation is necessary to assess the inlet error correctly; coarse grids can mask the boundary-layer defect and produce artificially optimistic $H^1$-estimates. Under noise, KLO reconstructions remain smoother and structurally consistent in the interior, but require careful tuning of the regularisation prefactor and smoothing width.

The ensemble statistical study ($n=1000$ realisations, SG-to-KLO transfer setting) reveals a clear noise-dependent crossover. At $\delta=0\%$, SG outperforms KLO in both metrics across all realisations (Cohen's $d=-1.54$ for $H^1$, $-3.26$ for $L^2$; KLO win rate $0\%$). At $\delta=1\%$, neither method dominates in $H^1$: the Wilcoxon signed-rank test does not reject the null hypothesis ($p=0.119$), the rank-biserial correlation is $-0.057$, and KLO attains the smaller error in $55.8\%$ of realisations, while KLO gains a moderate $L^2$ advantage (win rate $67.6\%$). At $\delta=5\%$ and $10\%$, KLO is clearly superior in $H^1$ (Cohen's $d>2.2$, win rate $>99\%$) and maintains a consistent $L^2$ advantage (win rate $85$--$90\%$). This crossover reflects a fundamental trade-off: the method with the highest noiseless accuracy is not necessarily the most noise-robust.

In summary, SG is the recommended choice when data are clean and computational simplicity matters. KLO is preferable when measurements are noisy ($\delta\gtrsim1\%$) and robustness to perturbations outweighs achieving the smallest noiseless error, provided inlet regularisation and grid resolution are carefully controlled.

\paragraph{Control of numerical and statistical bias.}
Several precautions were taken to ensure a fair comparison between SG and KLO. For each realisation, the same area profile was used for both methods, and the forward data were computed on a refined grid before being resampled onto a common measurement grid. Both methods were supplied with the same SG-generated boundary trace at each noise level, with identical additive noise realisations. The true profile and both reconstructions were then evaluated on a common spatial grid before computing the relative \(L^2\)- and \(H^1\)-errors; in particular, the KLO reconstruction was interpolated onto the SG reconstruction grid. Thus, the statistical tests were performed on paired errors obtained from matched profiles, data, noise levels, and discrete error norms.

The comparison should therefore be interpreted as an assessment of the implemented regularised reconstruction pipelines under matched numerical conditions. In particular, KLO uses smoothing and admissible range constraints to stabilise the differentiation of the reconstructed control quantity, while SG uses its own Fredholm regularisation. 

\section{Acknowledgements}
This work was supported by the Finnish Ministry of Education and Culture’s Pilot for Doctoral Programmes (Pilot project Mathematics of Sensing, Imaging and Modelling) and by the Research Council of Finland through the Flagship of Advanced Mathematics for Sensing, Imaging and Modelling (decision number 359183).

\printbibliography

\end{document}